\documentclass[11pt]{article}
\author{Bing-Long Chen\\[8pt]}
\title{\textbf{On stationary solutions to the vacuum Einstein field equations}}
\date{April 4, 2016}
\usepackage{latexsym}
\usepackage{amsmath}
\usepackage{amssymb,amscd}
\usepackage[mathscr]{eucal}
\newtheorem{thm}{Theorem}[section]
\newtheorem{cor}[thm]{Corollary}
\newtheorem{lem}[thm]{Lemma}

\newtheorem{rem}{Remark}[section]

\numberwithin{equation}{section}
\newenvironment{pf}{{\noindent \it  Proof.}}{{\hfill$\Box$}\\}
\begin{document}
\maketitle
\let\thefootnote\relax\footnotetext{AMS Mathematics Subject Classification Numbers: Primary
53c50; Secondary 83c20. }

\begin{abstract}
We prove that  any  4-dimensional geodesically complete spacetime with a timelike Killing field satisfying  the vacuum  Einstein field equation  $Ric(g_{M})=\lambda g_{M}$ with nonnegative cosmological constant $\lambda\geq 0$  is flat. When dim $\geq 5$, if  the spacetime is  assumed to be static additionally,  we  prove that its universal cover  splits isometrically as  a product of a Ricci flat Riemannian manifold and  a real line. 
\end{abstract}

  \section{Introduction}

   A Lorentzian manifold  $(M,g_{M})$ or a spacetime is a differentiable manifold $M$ equipped with a Lorentzian metric $g_{M}$  of signature $(-1,+1,\cdots, +1)$.    In general relativity,  the gravity is described by a spacetime  4-manifold $(M,g_{M})$,  the Lorentzian metric $g_{M}$ satisfies  the Einstein field equation:
   \begin{equation}\label{EE}
   \begin{split}
   Ric(g_{M})-\frac{1}{2}R g_{M}+\varLambda g_{M}=\kappa T
   \end{split}
   \end{equation}
   where $T$ is the energy-momentum tensor due to the presence of matter or fields, $\kappa$ and $\varLambda$ are constants.   
   
   In this paper, we are interested in the solutions to (\ref{EE}) with  timelike Killing fields.  These solutions are called stationary solutions. Stationary solutions  are  used to  model the possible time-independent limit states of a  cosmological  system.  For instance, Kerr metrics  are   stationary and vacuum solutions ($T=0,\Lambda=0$) to (\ref{EE}), while Schwarzschild metrics are static and vacuum solutions ($T=0,\Lambda=0$). Here, static means that the spacetime has a timelike Killing field whose  orthogonal complement is an integrable distribution, i.e., the timelike Killing field is locally orthogonal to spacelike hypersurfaces.  These stationary solutions, including Schwarzchild, Kerr, Reissner-Nordstrom (electrovac static), Kerr-Newmann metrics (electrovac stationary),    have been central to the study of the black hole spacetimes, see \cite{Chr} \cite{Ro}. 
   
   If  the spacetime   $(M,g_{M})$   admits an isometric $R$-action such that the $R$-orbits are timelike curves, the infinitesimal generator of the $R$-action  is  a timelike Killing field.  In many literatures,  the terminology " stationary"  was also used to referring to  the existence of such global $R$-action. Here,  our usage of  "stationary"  is in a broader sense, it only refers to the existence of a timelike Killing field.

  One of the main results of this  paper is  the following theorem:  
 \begin{thm} \label{t1.1} Let $(M,g_M)$ be a geodesically complete spacetime of dimension 4  with a timelike Killing field $X$ such that  $g_{M}$ satisfies the  Einstein equation $Ric(g_M)=\lambda g_{M}$, where $\lambda\geq 0$.  Then $(M,g_M)$ is  flat. 
 \end{thm}
 
  Here,  $(M,g_M)$ is said to be geodesically complete if the affine parameters of  any $g_M$-geodesic on $M$ can be extended to the whole real line $R$ .  The Einstein equation satisfied by the spacetime in  Theorem \ref{t1.1}  is   equivalent to $T=0$ and $\Lambda\geq 0$ in (\ref{EE}). We remark that when $\lambda<0$, the result of Theorem \ref{t1.1} is not true.  The simplest counter examples are  anti-De Sitter spacetimes, which are static, geodesically complete, and satisfying $Ric(g_{M})=\lambda g_{M}$ for $\lambda<0$.

    Recall that a Lorentzian manifold $(M,g_{M})$ is said to be chronological if it contains no closed timelike curves. In \cite{A}, M. T.  Anderson proved that if  the spacetime $(M^4,g_{M})$ is geodesically complete, chronological and admits an  isometric timelike $R$-action such that $g_{M}$ satisfies the vacuum Einstein field equation $Ric(g_{M})\equiv 0$,  then  $(M^4,g_{M})$ must be flat.    When the $R$-orbit space $M/R$ is an asymptotically flat 3-manifold, the result was due to A.  Lichnerowicz \cite{L} in 1955. The  previous pioneering work was due to A. Einstein and A. Einstein-W. Pauli, see \cite{EP}.

     The  asymptotic flatness on the orbit space is  usually a reasonable assumption for an isolated chronological  physical system.  
  The chronological condition is  used to ensure that the $R$-orbit space $M/R$ (denoted by $N$) is a   paracompact  Hausdorff  and smooth manifold, see \cite{H}. Actually, in this case,  the manifold $M$ is diffeomorphic to $R\times N$, and the metric $g_M$ has the following global form (see \cite{A},\cite{H}, \cite{K}) 
  \begin{equation}\label{met}\begin{split}
  g_{M}=-u^2(dt+\pi^{\ast}\theta)^2+ \pi^{\ast}g_{N}, 
\end{split}
\end{equation} 
on $M\approx R\times N$, where $u$, $\theta$ are some function and $1-$form on $N$, $g_{N}$ is a Riemannian metric on $N$, $\pi: M\rightarrow N$  is the projection map from $M$ to the space of $R$-orbits $N$. The   argument in \cite{A}  used the collapsing theory(c.f.\cite{CG3}) for  a sequence of 3-Riemannian manifolds which are the orbit spaces  of the isometric $R-$actions.  When the orbit spaces $N$ are noncompact and have dimension equal to 3,  Anderson \cite{A}  argued that the collapsing can be unwrapped by considering their  universal covers.  Recently, J. Cortier and V. Minerbe \cite{CM} gave a new proof of Anderson's theorem \cite{A} under an extra assumption on the norm of the timelike Killing field $X$.  

 Without chronological condition, the orbit space could be very "bad".  A simple compact example is Minkowski  flat torus $T^2$ (see \cite{H}), here we  take the constant vector field with irrational slope as a timelike Killing field.  In this case, any Killing  orbit is dense in $T^2$, so the  quotient topology just consists of two elements: the empty set and the whole space.   For noncompact examples, one can take a product of such a torus with a real line.

Whether the chronological condition can be removed in  Anderson's theorem is a question asked in \cite{A} (see \cite{A} \S1 second paragraph), so Theorem \ref{t1.1} answers this question affirmatively.

 Direct generalization of Theorem \ref{t1.1} to higher dimensions is not true, because we have to allow  non-flat examples which are  product of a Ricci flat Riemannian manifold  with  a real line.  So when dimension $\geq 5$, the best we can hope  is a splitting result.  Actually,  if   the  spacetimes are  assumed to be static,  we can prove that it is really the case:  
  \begin{thm} \label{t1.2} Let $(M,g_M)$ be a geodesically complete spacetime of dimension $n+1$ with a timelike Killing field whose orthogonal complement is  integrable.  Suppose the metric  $g_{M}$ satisfies  the   Einstein  equation $Ric(g_M)=\lambda g_{M}$, where $\lambda\geq 0$. Then $\lambda=0$ and the universal cover of $(M,g_M)$ is isometric to $R\times N$ equipped with a product metric $-dt^2+g_{N}$, where $(N,g_{N})$ is a complete Ricci flat Riemannian manifold of dimension $n$. 
 \end{thm}

  As we mentioned before, the result in Theorem \ref{t1.2} is not true for $\lambda<0$.  
  
  It should be noted that recently M. Reiris \cite{RM} has shown Theorem \ref{t1.2}  under the chronological condition.  More precisely,   M.  Reiris  \cite{RM} has obtained the same result  for  static solutions to Einstein-scalar equation under the assumption that the spacetime splits  topologically  as $M\approx R\times N$ and the metric has global form (\ref{met}) with $\theta=0$.   
  
   Theorems \ref{t1.1} and \ref{t1.2} are  derived by proving a  local curvature estimate or a local  gradient  estimate of the norm of the  Killing field $X$. To state the result, we need to introduce a Riemannian metric  which is naturally associated to the stationary spacetime $(M,g_{M},X)$.

   Let $X^{\ast}$ be the  1-form on $M$ obtained from $X$ by lowering indices.   We define  \begin{equation}\label{hat}
\hat{g}=-\frac{2}{g_{M}(X,X)}X^{\ast} \otimes X^{\ast}+g_{M}, 
\end{equation} 
 which is a Riemannian metric on $M$.  It can be shown that the vector field $X$ is still a Killing field for the metric  $\hat{g}$. In other words, one can associate a   stationary Riemannian metric  $\hat{g}$ to a stationary Lorentzian metric   $g_{M}$ with the same Killing field.  See  \cite{CL1} and  \cite{CL2} for similar ideas  in treating the injectivity radius estimate  and local optimal regularity of Einstein spacetimes.

 Our local curvature or gradient estimates are  the followings:  
 \begin{thm}\label{t1.3}Let $(M,g_M)$ be a   spacetime  of dimension $4$ with a timelike Killing field $X$ and $g_{M}$ satisfies  the  Einstein  equation $Ric(g_M)=\lambda g_{M}$.    Let $\hat{B}(x_0, a)$ be a $\hat{g}-$metric ball centered at $x_0$ of radius $a>0$ with compact  closure  in $M$. Then there is a universal constant $C>0$  such that 
 \begin{equation} \label{ce}
\begin{split}
\sup_{x\in \hat{B}(x_0,\frac{a}{2})}\mid Rm(g_M)\mid_{\hat{g}}\leq C(a^{-2}+\max\{-\lambda,0\}).
\end{split}
\end{equation}
 \end{thm}

 \begin{thm}\label{t1.4}Let $(M,g_M)$ be a  spacetime of dimension $n+1$ with a timelike Killing field $X$ whose orthogonal complement is integrable, and $g_{M}$ satisfies  Einstein  equation $Ric(g_M)=\lambda g_{M}$.   Let $\hat{B}(x_0, a)$ be a $\hat{g}-$metric ball centered at $x_0$ of radius $a>0$ with compact closure in  $M$. Then there is a universal constant $C>0$  such that 
 \begin{equation} \label{ge}
\begin{split}
& \sup_{x\in \hat{B}(x_0,\frac{a}{2})}\mid\hat{\nabla}\log (-g_{M}(X,X))\mid_{\hat{g}}\leq C(\sqrt{n}a^{-1}+\sqrt{\max\{-\lambda,0\}}).
\end{split}
\end{equation}
 \end{thm}
Note that  $\max\{-\lambda,0\}=0$  if $\lambda\geq 0$  in  Theorems \ref{t1.3} and \ref{t1.4}. 
 To prove Theorems \ref{t1.1} and \ref{t1.2} from Theorems \ref{t1.3} and \ref{t1.4}, we need a fact that $g_{M}$- geodesic  completeness implies $\hat{g}$-geodesic completeness (see Theorem \ref{complete} in Section 3).  
 
 When dimension equals to $4$, actually we  can show that a local curvature estimate holds on more general spacetimes which are not necessarily vacuum (see Theorems \ref{N},  \ref{M}). The result roughly says that if the energy momentum tensor is controlled, then the full curvature tensor of the spacetime can also be controlled quantitatively. The non-vacuum Einstein field equation coupled with specific matter fields will be treated in forthcoming papers.      
 
 The paper is organized as follows. In section 2, we prepare some preliminary formulas that will be used throughout the paper.  Sections 2.1, 2.2 and 2.3 involve many straight forward computations on the connections and curvatures, and the formulas work for stationary Lorentzian and stationary  Riemannian manifolds.  In section 3, we prove that the $g_M$-geodesic completeness implies the  $\hat{g}$-geodesic completeness. In section 4, we prove  Theorems \ref{t1.4} and \ref{t1.2}. Theorems  \ref{t1.3} and  \ref{t1.1} are proved in section 5.   
 
 \textbf{Acknowledgement } The author is grateful to Professors  S. T. Yau, X. P. Zhu and Dr. J. B. Li for helpful discussions. Professor S. T. Yau told the author  that  the estimate (at least $\lambda=0$ case)  in Corollary 4.3 was also known to Professor R. Schoen. 
 The work was partially supported  by grants NSFC11521101, 11025107. 

\section{Stationary spacetime and its  associated  Riemannian metric}

Suppose $(M, g_M,X)$  is a stationary spacetime of dimension $n+1$, where $g_{M}$ is a smooth  Lorentzian metric on $M$ and $X$ is a timelike Killing field. Denote the set  of  integral curves of $X$ by $N$,  $\pi: M\rightarrow N$ the projection map.  

\subsection{A local coordinate system}
Fix a point $P\in M$,  we will construct a natural coordinate system $\{x^{\alpha}\}$ around $P$ in the followings so that the metric has the form (\ref{met}) locally.

Let  $\Psi_{\tau}$  be   the  (local) diffeomorphisms generated by $X$ such that $\Psi_{0}=id$, and $\Psi_{\tau_1}\Psi_{\tau_2}=\Psi_{\tau_1+\tau_2}$ wherever they are defined. 
Now we fix  a codimensional one spacelike submanifold $\Sigma\subset M$  passing through $P$ such that $\bar{\Sigma}$ is compact.  Considering  the affine parameters  of all integral curves  of $X$ starting from $\Sigma$, we obtain   a  function $t$ defined on an open neighborhood of  $\Sigma$ in  $M$  such that $t=0$ on $\Sigma$.  Given a local coordinate system $(x^1,\cdots, x^n)$ on  $\Sigma$ around $P$.  We can construct a local coordinate system $(x^0,x^1,\cdots,x^n)$  on $M$ around $P$ , where $x^0=t$.  Actually, for any point $Q\in \Sigma$ lying in the coordinate chart in $\Sigma$, we require  $x^{i}(\Psi_t(Q))=x^{i}(Q)$ for $i=1,2\cdots,n$, and $x^0(\Psi_{t}(Q))=t$.  Throughout the paper, we  use  Greek  letters $\alpha, \beta, \cdots$  to indicate the indices varying  from $0$ to $n$ and Latin letters $i,j,k,\cdots$ varying from $1$ to $n$.  The coordinate system $\{x^{\alpha}\}$ depends on the choice of spacelike submanifold $\Sigma$ and the coordinate system $\{x^{i}\}$ on $\Sigma$.

    Let $X^{\ast}$ be the 1-form obtained by lowering indices of  $X$. The induced Riemannian metric on the horizontal distribution $\mathcal{H}=X^{\perp}$ is given by  $g_{\mathcal{H}}=g_{M}-\frac{1}{g_{M}(X,X)}X^{\ast} \otimes X^{\ast}$. It is clear that  the horizontal metric $g_{\mathcal{H}}$ and  1-form  $\theta=-u^{-2}(X^{\ast}+u^{2}dt)$ satisfy  $\mathcal{L}_{X} g_{\mathcal{H}}=0$, $\mathcal{L}_{X}\theta=0$,  where $u^2=-g_{M}(X,X)$,  $\mathcal{L}_{X}$ is the Lie derivative of the vector field  $X$.  
    
    It we choose a different spacelike submanifold, say $\Sigma'$,  and denote the corresponding time function and one-form   by $t'$ and $\theta'$, we have $\theta-\theta'=d \psi$, where $\psi=t'-t$ is a locally defined smooth  function.  Note that  the integrability of the horizontal distribution $\mathcal{H}$ is equivalent to $dX^{\ast}=0$ mod $X^{\ast}$ (Frobenius condition), or $d\theta=0$.

         The metric $g_M$ now has the following form  
  \begin{equation}\label{2.1}\begin{split}
  g_{M}=-u^2(dt+\theta)^2+ g_{\mathcal{H}}, 
\end{split}
\end{equation} 
on a neighborhood of $P$.  In the above local coordinate system $\{x^{\alpha}\}$,  we have 
\begin{equation}\label{2.1'}\begin{split}
&g_{\mathcal{H}}(\frac{\partial}{\partial x^0},\frac{\partial}{\partial x^{\alpha}})=\theta(\frac{\partial}{\partial x^0})=0,  \\
& \frac{\partial}{\partial t} g_{ij}= \frac{\partial}{\partial t}\theta_{i}= \frac{\partial}{\partial t}u^2=0,
\end{split}
\end{equation} 
where  $g_{ij}\triangleq g_{\mathcal{H}}(\frac{\partial}{\partial x^{i}},\frac{\partial}{\partial x^{j}})$, $\theta_i=\theta(\frac{\partial}{\partial x^i})$.  
Roughly speaking, the  equations in   (\ref{2.1'}) say  that  $u$, $\theta$,  $g_{ij}$  are essentially quantities on the space $N$ of $X$-integral curves. Actually, if we 
identify $\Sigma$ with  $\pi(\Sigma)$ using the projection map $\pi: M\rightarrow N$ and  equip $\pi(\Sigma)\subset N$  the Riemannian structure from  $(\Sigma,g_{\mathcal{H}})$. The equations in  (\ref{2.1'})  are  equivalent to   $u=\pi^{\ast} u$, $\theta=\pi^{\ast}(\theta)$, $\pi^{\ast}g=g$.

\subsection{Connection and curvature matrices}
Now we will do  some straight forward  calculations for  metrics  of  the form $\bar{g}=w(dt+\theta)^2+g$, where $w,\theta,g$ are $t-$independent.  The metric is Lorentzian if $w<0$,  and  Riemannian if $w>0$.

 Since $\bar{g}_{00}=w$, 
 $\bar{g}_{0i}=w\theta_i$, $\bar{g}_{ij}=g_{ij}+w\theta_i\theta_j$, one  can calculate  the inverse matrix $(\bar{g}^{\alpha\beta})$ of $(\bar{g}_{\alpha\beta})$: 
 \begin{equation}\label{inverse}
 \begin{split}
 &\bar{g}^{00}={w}^{-1}+|\theta|^2, \bar{g}^{0i}=-\theta^{i}, \bar{g}^{ij}=g^{ij},
 \end{split}
 \end{equation} where $\theta^i=g^{ij}\theta_j$ and $|\theta|^2=g^{ij}\theta_i\theta_j$.  
 It is useful to choose a good frame to calculate the connection coefficients.  Let $e_0=\frac{\partial}{\partial t}$,  $ e_i= \frac{\partial}{\partial x^i}-\theta_i \frac{\partial}{\partial t}$.  It can be shown  that  $[e_0,e_j]=0$, $[e_i,e_j]=-\Lambda_{ij}e_0$ and $\langle e_0, e_i\rangle=0$,  $\langle e_i,e_j\rangle=g_{ij}$, where  $\Lambda_{ij}=\nabla_{i}\theta_j-\nabla_{j}\theta_i$, $\nabla_{i}\theta_j$ is the covariate derivative of the tensor $\theta$ w.r.t. the horizontal metric $g$.

  The dual frame of $\{e_{\alpha}\}$ is $\{\omega^{\alpha}\}$, $\omega^{0}=dt+\theta, \omega^{i}=dx^i, i=1,2,\cdots,n$.  They satisfy $\langle \omega^0, \omega^0\rangle=w^{-1}$,  $\langle \omega^0, \omega^i\rangle=0$,   $\langle \omega^i, \omega^j\rangle=g^{ij}$.    Denote the Levi-Civita connection matrix w.r.t. the basis  $e_{\alpha}$ by ${\bar{\omega}}^{\beta}_{\alpha}$, $De_{\alpha}=\bar{\omega}^{\beta}_{\alpha}\otimes e_{\beta}$.  Recall  Cartan's   equations, 
  \begin{equation}\label{cfm}
  \begin{split}
  d\omega^{\alpha}&=\omega^{\beta}\wedge\bar{\omega}^{\alpha}_{\beta}\\
  d\langle \omega^{\alpha},\omega^{\beta} \rangle&=-\bar{\omega}_{\gamma}^{\alpha}\langle \omega^{\gamma},\omega^{\beta} \rangle-\bar{\omega}_{\gamma}^{\beta} \langle \omega^{\gamma},\omega^{\alpha} \rangle\\
  \bar{\Omega}_{\alpha}^{\beta}&=d\bar{\omega}_{\alpha}^{\beta}-\bar{\omega}_{\alpha}^{\gamma}\wedge\bar{\omega}_{\gamma}^{\beta}.
  \end{split}
  \end{equation}  
   The 1st equation  in (\ref{cfm}) says that the connection is torsion free, the 2nd  says it is compatible with  the metric $\bar{g}$.   The 3rd equation in (\ref{cfm}) is the definition of the curvature matrix $\{\bar{\Omega}^{\beta}_{\alpha}\}$. 
   The connection matrix $\{\bar{\omega}_{\alpha}^{\beta}\}$ is completely determined by the first two equations in (\ref{cfm}). Actually,  the 2nd  equation (\ref{cfm}) takes the following form 
    \begin{equation}
   \begin{split}\label{cpme}
   dg^{ij}=-\bar{\omega}^i_{k}g^{kj}-\bar{\omega}^i_{k}g^{kj}\\
   0=-\bar{\omega}^{0}_kg^{ki}- \bar{\omega}^{i}_{0}w^{-1}\\
   dw^{-1}=-2\bar{\omega}^{0}_{0}w^{-1}.
   \end{split}
   \end{equation}
   The 1st  equation in (\ref{cfm}) is 
  \begin{equation}
   \begin{split}\label{tf}
   d\theta=d\omega^0=\omega^k\wedge \bar{\omega}^{0}_k+\omega^{0}\wedge \bar{\omega}^{0}_{0}\\
   0=d\omega^i=\omega^k\wedge \bar{\omega}^{i}_k+\omega^{0}\wedge \bar{\omega}^{i}_{0}.
   \end{split}
   \end{equation}
Combining (\ref{cpme}) and (\ref{tf}),   we have 
   \begin{equation}\label{cfm1}
  \begin{split}
  \bar{\omega}_{j}^{i}&=\omega^{i}_{j}+\frac{1}{2}wg^{il}\Lambda_{jl}\omega^0\\
 \bar{\omega}_{0}^{i}&=-\frac{1}{2}wg^{il}\Lambda_{lk}\omega^k-\frac{1}{2}g^{il}\nabla_{l}w\omega^0\\
   \bar{\omega}_{i}^{0}&=\frac{1}{2}\Lambda_{il}\omega^l+\frac{1}{2}w^{-1}\nabla_{i}w\omega^0\\
    \bar{\omega}_{0}^{0}&=\frac{1}{2}w^{-1}dw,  \end{split}
  \end{equation}  
  where $\{\omega^{i}_{j}\}$   is the connection matrix of the horizontal metric $g=g_{\mathcal{H}}$  w.r.t. the natural frame $\{\frac{\partial}{\partial x^i}\}$. Now one  can calculate the curvature matrix by using (\ref{cfm1}) and the 3rd equation in (\ref{cfm}):
   \begin{equation}\label{cfm5}
  \begin{split}
  \bar{\Omega}_{j}^{i}&=\Omega^{i}_{j}+\frac{w}{4}g^{il}(\Lambda_{jp}\Lambda_{lq}+\Lambda_{jl}\Lambda_{pq})\omega^p\wedge\omega^q\\&\ \ +[\frac{1}{2}g^{il}\Lambda_{jl}dw+\frac{1}{2}wg^{il}\nabla_p\Lambda_{jl}\omega^p+\frac{1}{4}(\Lambda_{jl}\omega^l\nabla^iw-\nabla_jw\Lambda_{lp}g^{il}\omega^p)]\wedge \omega^0\\
 \bar{\Omega}_{0}^{i}&=g^{il}[-\frac{w}{4}\nabla_l\Lambda_{pq}-\frac{1}{8}(w_p\Lambda_{lq}-w_q\Lambda_{lp})-\frac{1}{4}\nabla_lw\Lambda_{pq}]\omega^p\wedge\omega^q\\
 &\ \ +g^{il}(-\frac{1}{2}\nabla_{pl}w+\frac{1}{4}w^{-1}\nabla_pw\nabla_lw+\frac{w^2}{4}\Lambda_{pq}\Lambda_{lm}g^{qm})\omega^p\wedge \omega^0.\\
  \end{split}
  \end{equation} 
     Using the formula $De_{\alpha}=\bar{\omega}_{\alpha}^{\beta}\otimes e_{\beta}$,  (\ref{cfm1}) may be paraphrased  as follows: 
  
   \begin{equation} \label{cov}
\begin{split} 
& D_{e_i}e_j=\Gamma^{k}_{ij}e_k-\frac{1}{2}\Lambda_{ij}e_0\\
& D_{e_0}e_i=D_{e_i}e_0=\frac{1}{2}w\Lambda_{ik}g^{kl}e_l+\frac{1}{2}\nabla_{i}\log |w| e_0\\
& D_{e_0}e_0=-\frac{1}{2}g^{ij}\nabla_{i}we_j.\\
\end{split}
\end{equation}    
 
  Using  $\Omega^{\alpha}_{\beta}=\frac{1}{2}\bar{R}^{\alpha}_{\beta\gamma\delta}\omega^{\gamma}\wedge \omega^{\delta}$ and $\bar{R}(e_{\alpha},e_{\beta},e_{\gamma},e_{\delta})=\langle e_{\alpha},e_{\epsilon}\rangle\bar{R}^{\epsilon}_{\beta\gamma\delta}$, (\ref{cfm5}) can be rewritten as 
\begin{equation} \label{curvature}
\begin{split}
& \bar{R}(e_i,e_j,e_k,e_l)=R_{ijkl}+\frac{w}{4}(\Lambda_{il}\Lambda_{jk}-\Lambda_{ik}\Lambda_{jl})-\frac{w}{2}\Lambda_{ij}\Lambda_{kl}\\& \bar{R}(e_i,e_j,e_k,e_0)=-\frac{1}{2}(w\nabla_{k}\Lambda_{ij}+{\nabla_{k}w}\Lambda_{ij}) +\frac{1}{4}(\nabla_{i}w \Lambda_{jk}-\nabla_{j}w \Lambda_{ik})\\
& \bar{R}(e_i,e_0,e_j,e_0)=-\frac{1}{2}\nabla_{ij}w+\frac{1}{4}w^{-1}\nabla_{i}w\nabla_j w+\frac{w^2}{4}\Lambda_{ik}\Lambda_{jl}g^{kl}.
\end{split}
\end{equation}
  Here, our convention for the sign of the curvature tensor  is that we require $R_{ijij}>0$ on spheres.  (\ref{curvature}) can also be obtained alternatively by using (\ref{cov}) and the following formula: $$\bar{R}(e_{\alpha},e_{\beta},e_{\gamma},e_{\delta})=-\langle (D_{e_{\alpha}}D_{e_{\beta}}-D_{e_{\beta}}D_{e_{\alpha}}-D_{[e_{\alpha},e_{\beta}]})e_{\gamma},e_{\delta}\rangle.$$
   

Since we have computed all connection coefficients (see (\ref{cov})),   it is not difficult to compute the Hessian and the Laplacian of any time-independent function $f$: 
 \begin{equation} \label{covf}
\begin{split}
& \bar{\nabla}^2f(e_0,e_0)=\frac{1}{2}\langle\nabla w,\nabla f\rangle \\
&\bar{\nabla}^2f(e_0,e_j)=-\frac{1}{2}wg^{kl}\Lambda_{jk}f_{l}\\
& \bar{\nabla}^2f(e_i,e_j)=\nabla_{ij}f
\end{split}
\end{equation}
and 
 \begin{equation} \label{lf}
\begin{split}
& \bar{\triangle}f=\triangle f+\frac{1}{2}\langle \nabla \log|w|,\nabla f\rangle.
\end{split}
\end{equation}
The formulas (\ref{covf}) and (\ref{lf}) are  important in the calculations  of sections 4 (see (\ref{3.13})  and 5. 
\subsection{Ricci curvature}
 
 By taking traces on (\ref{curvature}), we get the Ricci curvature formula:  
  \begin{equation} \label{Ric}
\begin{split}
& \bar{R}ic(e_0,e_0)=-\frac{\triangle w}{2}+\frac{|\nabla w|^2}{4w}+\frac{w^2}{4}| \Lambda|^2 \\
&\bar{R}ic{(e_0,e_j)}=\frac{w}{2}g^{kl}(\nabla_k\Lambda_{jl}+\frac{3}{2}\Lambda_{jk}\nabla_{l}\log w) 
\\& \bar{R}ic({e_i,e_j})=R_{ij}-\frac{\nabla_{i}\nabla_j w}{2w}+\frac{\nabla_i w\nabla_jw}{4w^2}-\frac{w}{2}g^{kl}\Lambda_{ik}\Lambda_{jl},
\end{split}
\end{equation}
where $R_{ij}$ is the Ricci curvature of  the horizontal metric $g_{ij}$.

Let $w=-u^2<0$ in (\ref{2.1}) and (\ref{Ric}),  we have  
\begin{equation} \label{equ}
\begin{split}
& \triangle u =-\frac{u^3}{4}\mid \Lambda\mid^2+ u^{-1}\bar{R}ic(e_0,e_0) \\
& g^{kl}(\nabla_k\Lambda_{jl}+3\Lambda_{jk}\nabla_{l}\log u) =-2u^{-2}\bar{R}ic(e_0,e_j)
\\& R_{ij}=u^{-1}{\nabla_{i}\nabla_j u}-\frac{u^2}{2}g^{kl}\Lambda_{ik}\Lambda_{jl}+\bar{R}ic(e_i,e_j).
\end{split}
\end{equation}
Let \begin{equation}\label{hat1}
\hat{g}=-\frac{2}{g_{M}(X,X)}X^{\ast} \otimes X^{\ast}+g_{M}
\end{equation} 
be the  Riemannian metric defined in (\ref{hat}) on $M$. Combining (\ref{equ}) and (\ref{Ric}), we get
 \begin{equation} \label{hatric}
\begin{split}
& \hat{R}ic({e_0,e_0})=\frac{u^4}{2}| \Lambda|^2-\bar{R}ic(e_0,e_0) \\
& \hat{R}ic({e_0,e_j})=-\bar{R}ic(e_0,e_j)
\\& \hat{R}ic({e_i,e_j})=-u^2g^{kl}\Lambda_{ik}\Lambda_{jl}+\bar{R}ic(e_i,e_j),
\end{split}
\end{equation}
where $\hat{R}ic$ and $\bar{R}ic$ are Ricci curvatures of metrics $\hat{g}$ and $\bar{g}$.

 It is  helpful  to introduce a  new metric $\tilde{g}$ conformal to $g$ on horizontal distribution. This metric will play an important role in a priori estimates(see Section 5.3).   Let $\tilde{g}=u^{\frac{2}{n-2}}{g}$ be a conformal change of the horizontal metic ${g}$.  The Christoffel symbols of $\tilde{g}$ can be given by (see Chapter 5 in  \cite{SY}) 
 \begin{equation} \label{Chtil}
\begin{split}
 \tilde{\Gamma}^{k}_{ij}=\Gamma^{k}_{ij}+\frac{1}{n-2}(\nabla_i \log u\delta_{j}^{k}+\nabla_{j}\log u\delta^{k}_{i}-g^{kl}\nabla_{l}\log u g_{ij}).
\end{split}
\end{equation}
This implies 
 \begin{equation} \label{laplacetil}
\begin{split}
 u^{\frac{2}{n-2}}\tilde{\triangle} L=\triangle L+\langle \nabla \log u,\nabla L\rangle=\hat{\triangle}L
\end{split}
\end{equation}
for any t-independent smooth function $L$ on $M$.

The Ricci curvature of $\tilde{g}$ can be computed by the following formula (see Chapter 5 in  \cite{SY}): 
 \begin{equation} \label{Rictil}
\begin{split}
& \tilde{R}_{ij}={R}_{ij}-{\nabla}^2_{ij}\log  u+\frac{1}{n-2}(\log u)_{i}(\log u)_{j}
 -\frac{1}{n-2}\frac{\triangle u}{u}{g}_{ij}\\
 &\ \ \ \ \  =\frac{u^2}{4(n-2)}|\Lambda|^2 g_{ij}-\frac{u^2}{2} g^{kl}\Lambda_{ik}\Lambda_{jl}+\frac{n-1}{n-2}\frac{u_iu_j}{u^2}\\ & \ \ \ \  \ \ \ +\bar{R}ic(e_i,e_j)-\frac{u^{-2}\bar{R}ic(X,X)}{n-2}g_{ij},
 \end{split}
\end{equation} 
where we have used (\ref{equ}). 
By direct computations, we also  have 
\begin{equation} \label{equtil}
\begin{split}
& u^{\frac{2}{n-2}}\tilde{\triangle}\log  u =-\frac{u^2}{4}\mid \Lambda\mid^2+ u^{-2}\bar{R}ic(X,X) \\
& g^{kl}(\tilde{\nabla}_k\Lambda_{jl}+(3+\frac{4-n}{n-2})\Lambda_{jk}\nabla_{l}\log u) =-2u^{-2}\bar{R}ic(e_0,e_j).
\end{split}
\end{equation}

\begin{cor} \label{c3.1}
Let $(M,g_{M},X)$ be  a  static Einstein spacetime of dimension $n+1$ satisfying  $Ric(g_{M})=\lambda g_{M}$.  Then $(M,\hat{g})$ is a Riemannian  Einstein manifold with the same cosmological constant as $g_{M}$, i.e.,  
$Ric(\hat{g})=\lambda \hat{g}$. 
\end{cor}
\begin{pf} This follows from (\ref{hatric}) by noting that $d\theta=0$ on static spacetimes. 
\end{pf}

\begin{cor} \label{c3.2}
Let $(M,g_{M},X)$ be  a geodesically  complete stationary and chronological spacetime of dimension $n+1$ with a timelike Killing field $X$ such that  the  $X-$orbit space $N$ is a compact smooth manifold.  Then  the followings hold\\
i)   If ${R}ic_{M}(X,X)\leq 0$ holds everywhere,  then  $(M,g_{M},X)$ is static, $|X|^2_{g_{M}}\equiv const.$ and there is a closed 1-form $\theta$ on $N$ such that  $(M,g_{M})$ is isometric to a metric $-(d\tau+\theta)^2+g_{N}$ on $R\times N$. \\
ii) If $g_{M}$ is Einstein and static,  $Ric(g_{M})=\lambda g_{M}$,  we have  $\lambda=0$, the conclusion of i) holds,  and $(N,g_{N})$ is Ricci flat. 
\end{cor}
\begin{pf} 
First of all, by \cite{H} (c.f.\cite{A},\cite{K}),  $M$ is diffeomorphic to $R\times N$, and  the metric $g_{M}$ now has the  global form (\ref{met}). 
From the first equation of (\ref{equ}), we know $Ric_{M}(X,X)=u\triangle u+\frac{u^4}{4}|\Lambda|^2$. For  i), since $\triangle u\leq 0$, by strong maximum principle, we have  $Ric_{M}(X,X)=0$ and $d\theta=\Lambda=0$ and  $u\equiv const.$.  This shows i).  \\
 If we assume $Ric(g_{M})=\lambda g_{M}$, then $Ric_{M}(X,X)=-\lambda u^2$. By i), we know $\lambda\leq 0$. The first equation of (\ref{equ}) implies $\triangle u=-\lambda u$. Since $u>0$, we know $\lambda=0$ and $u=const.$ by strong maximum principle. The conclusion of ii) holds.  
 \end{pf}
\begin{rem} Under the assumptions of Corollary \ref{c3.2} and i),  $(M,g_{M})$ is isometric to a product  $-dt^2+g_{N}$ on $R\times N$ if $H^1(M,R)=0$.  
\end{rem}

\section{Completeness}
Before the discussion of completeness, we need to do some preliminary work on projecting   curves to horizontal ones.  Here, we say a curve is horizontal if its tangent vectors are horizontal.

\begin{lem} \label{proj} Let $(M,g_{M})$ be a spacetime with a timelike Killing field $X$. Let $\bar{\gamma}:I\rightarrow M$  be a  smooth curve on $M$, where $I\subset R$ is an interval. Fix $s_0\in I$, $p_0=\Psi_{\tau_0}(\bar{\gamma}(s_0))$, there is  a unique  maximal  smooth horizontal curve $\sigma:I'\rightarrow M$ such that $I'\subset I$, $\tau(s_0)=\tau_0$, $\Psi_{\tau(s)}(\bar{\gamma}(s))=\sigma(s)$ for any $s\in I'$, where $\tau:I'\rightarrow R$ is a  smooth function, $\Psi_{\tau}$ is the local flow generated by $X$.  Moreover, $I'=I$ provided that the vector field $X$ is complete. 
\end{lem}
Here $\sigma$ is maximal means that any such horizontal projection curve of $\bar{\gamma}$ passing $p_0$ is only a part of $\sigma$.  The vector field $X$ is said to be complete if   any integral curves of $X$ can be defined, for their affine parameters,  on  the whole real line $R$. 

\begin{pf}  Consider the curve  $\Psi_{\tau_0}(\bar{\gamma})$. Note that the map $\Psi_{\tau_{0}}$ might not be defined on whole $\bar{\gamma}$, so the parameters of the curve $\Psi_{\tau_0}(\bar{\gamma})$ lie in a connected subinterval of $I$.   Let $\Psi_{\tau_0}(\bar{\gamma}(s))$ be parameterized by $(x^0(s),x^1(s),\cdots,x^n(s))$  on a chronological chart $\{x^{\alpha}\}$ around $p_0$  used in  $\S{2.1}$.  Let  $T=\sum_{\alpha} T^{\alpha}e_{\alpha}$ be  the tangent vector of $\Psi_{\tau_0}(\bar{\gamma}(s))$, where $T^i=\frac{dx^{i}}{ds}$ and $T^0=\frac{dx^0}{ds}+\sum T^i\theta_i$. 
 
  Consider  another   curve \begin{equation}\label{2.9}\sigma(s)=(y(s),x^{1}(s),\cdots, x^{n}(s))\end{equation}  on the coordinate system $\{x^{\alpha}\}$, where we require $\frac{d}{ds}y(s)+\sum T^{i}\theta_i=0$, $\sigma(s_0)=\Psi_{\tau_0}(\bar{\gamma}(s_0))$. 
 It is clear 
 $ \frac{d \sigma(s)}{ds}=\sum_{i=1}^{n}T^ie_i$, i.e.,  $\sigma(s)$ is horizontal.  Note that the function $y(s)$ is determined uniquely by these requirements.  So $\bar{\gamma}$ has a unique horizontal projection $\Psi_{\tau(s)}(\bar{\gamma}(s))$ on this chart, where $\tau(s)=\tau_0+y(s)-x^0(s)$ is clearly a smooth function.  Because the manifold may be covered by chronological coordinate charts  used  in $\S{2.1}$, one  can extend the horizontal projection curve of $\bar{\gamma}$ to a maximal one. 
  \end{pf}
\begin{rem}  It is clear that  in Lemma \ref{proj},  we have \begin{equation} \label{bx}\begin{split}|\dot{\sigma}|_{\hat{g}}(s)\leq |\dot{\bar{\gamma}}|_{\hat{g}}(s),  \ \ s\in I', \end{split}\end{equation}
  which implies that if $\bar{\gamma}$ is not horizontal, then the $\hat{g}$-length of  the horizontal projection curve $\sigma$ will become strictly smaller. 
  \end{rem}
  
   \begin{lem} \label{Xcom}
Let $(M,g_{M})$ be  a time-geodesically complete spacetime  with a timelike Killing field $X$. Then $X$ is complete. 
\end{lem}
\begin{pf}  We only need to show that any integral curve  $\zeta:[a,b)\rightarrow M$ of $X$ can be extended over $b$.  When $c>a$ is close to $a$, there is a timelike geodesic $\gamma:[0,d]\rightarrow M$ such that $\gamma(0)=\zeta(a)$ and $\gamma(d)=\zeta(c)$. By time-completeness assumption, $\gamma$ can be extended to be defined on all affine parameters.  $\Psi_{b-c}$ is clearly defined near $\zeta(a)$ and $\Psi_{b-c}(\zeta(a))=\zeta(b-(c-a))$, where $\Psi_t$ are  the local diffeomorphisms generated by $X$.  To prove the lemma, it suffices to  show that the maps $\Psi_{\tau}$, for $\tau\in [0,b-c]$ can be defined on whole $\gamma$.  

 Suppose this is  not true, there will be a smooth family of $X-$integral curves connecting $\gamma(s)$ and $\sigma(s)$ for $s$ lying in a maximal interval $I=(a',b')$. Without loss of generality, we assume $b'<\infty$. The integral curve $\eta$ of $X$  starting at $\gamma(b')$ can only be defined on a maximal interval $[0, c')$ where $c'<b-c$. Considering the timelike geodesic $\xi(s)=\Psi_{c'}(\gamma)$ defined on $(a',b')$, it can also be extended to all affine parameters. Near the point $\xi(b')$, by considering  the integral curves of $-X$ starting at $\xi(s)$ for $s\in (b'-\epsilon,b']$, we find  $\eta$ can actually be extended  over $c'$, and $\eta(c')=\xi(b')$. This is a contradiction.  
 
\end{pf}

In general, the horizontal projection curve  of a geodesic is no longer a geodesic, but it still satisfies a "good" ODE.   Indeed, 
  in a coordinate system $\{x^{\alpha}\}$ in $\S{2.1}$, by using (\ref{cov}) and direct computations, we have 
  \begin{equation} \label{2.7}
\begin{split}
 \bar
 {\nabla}_{T}T= & w^{-1}\frac{d(T^0w)}{ds}e_0\\ &+[\frac{dT^{i}}{ds}+\Gamma^i_{kl}T^kT^l+\frac{1}{2}(T^0w)^2 \nabla_j w^{-1} g^{ij} +T^0w \Lambda_{lm}T^{l}g^{mi}]e_i,
\end{split}
\end{equation}
where $\bar{\gamma}(s)=(x^0(s),x^1(s),\cdots, x^{n}(s))$, $T=\dot{\bar{\gamma}}=\sum_{\alpha}T^{\alpha}e_{\alpha}$. 

 This implies that $\bar{\gamma}$ is a geodesic  on $(M,\bar{g})$ if and only if the curve $\gamma=\pi(\bar{\gamma})$, $\gamma(s)=(x^1(s),\cdots, x^{n}(s))$,  satisfies 
\begin{equation} \label{pge}
\begin{split}
 & T^0 w=\langle T,X\rangle=const.\triangleq c\\
 & \nabla_{\dot{\gamma}}\dot{\gamma}=-\frac{c^2}{2}\nabla w^{-1}-c(i_{\dot{\gamma}}d\theta)^{\sharp},
\end{split}
\end{equation}
where  $i_{\dot{\gamma}}d\theta$  is the   1-form  obtained from the contraction of 2-form $\Lambda=d\theta$  with the   tangent vector  $\dot{\gamma}=\sum T^i\frac{\partial }{\partial x^i}$.  

 The first equation of (\ref{pge}) can be derived alternatively  by $$T\langle T,X\rangle=\langle T,\bar{\nabla}_{T}X\rangle=0$$ 
 since $X$ is Killing and $\bar{\nabla}_{T}T=0$.

\begin{thm} \label{complete}
Let $(M,g_{M})$ be  a geodesically complete spacetime with a timelike Killing field $X$. Then $(M,\hat{g})$ is a complete Riemannian manifold, where $\hat{g}$ is defined in (\ref{hat}). 
\end{thm}
\begin{pf} Fix  $p\in M$, we will show that the exponential map $\exp_{\hat{g}}$ of $\hat{g}$ can be defined on  the whole tangent space $T_pM$  at  $p$. Let $\bar{r}$ be the supremum of all $r>0$ such that the exponential map $\exp_{\hat{g}}$ can be defined on a ball of radius $r$ centered at $0$ in $T_pM$. We have to show $\bar{r}=\infty$. 
We argue by contradiction.  Suppose $\bar{r}<\infty$.  Let    $\bar{\gamma}:[0,\bar{r})\rightarrow M$  be a normal $\hat{g}$-geodesic parameterized by arclength, $\bar{\gamma}(0)=p$, $|\dot{\bar{\gamma}}(0)|_{\hat{g}}=1$, and $\bar{\gamma}$ can not be extended  over time $\bar{r}$ (as a $\hat{g}$-geodesic).  From  the definition of $\bar{r}$, for any $r<\bar{r}$, we know $B(p,r)=\exp_{\hat{g}}(B(0,r))$ and $\bar{B}(p,r)=\exp_{\hat{g}}(\bar{B}(0,r))$, where $B(p,r)=\{q\in M: d_{\hat{g}}(q,p)<r\}$, $\bar{B}(p,r)=\{q\in M: d_{\hat{g}}(q,p)\leq r\}$,  $B(0,r)=\{v\in T_pM: |v|<r\}$,  and $\bar{B}(0,r)=\{v\in T_pM: |v|\leq r\}$.  Therefore,  $\bar{B}(p,r)=\{q\in M: d_{\hat{g}}(q,p)\leq r\}$ is compact provided $r<\bar{r}$.

 We assume that $\bar{\gamma}$ is not horizontal, otherwise $\bar{\gamma}$ can be extended to whole real line $R$, because a horizontal $\hat{g}$-geodesic is also  a horizontal $g_{M}$-geodesic (see (\ref{pge})). 
 By Lemmas \ref{proj},  \ref{Xcom},   one  can  construct a smooth horizontal projection curve $\sigma:[0,\bar{r})\rightarrow M$ such that  $\sigma(0)=\gamma(0)=p$, $\tau(0)=0$,  $\sigma(s)=\Psi_{\tau(s)}\gamma(s)$,  where  $\tau:[0,\bar{r})\rightarrow R$  is a smooth function.  
Since  $\bar{\gamma}$ is not horizontal,  the $\hat{g}$-length of $\sigma$ is less than that of $\gamma$ (see (\ref{bx})), i.e., $L\triangleq L(\sigma)<L(\gamma)=\bar{r}$.  By (\ref{bx}), for any $0<a<b<\bar{r}$, we have $d_{\hat{g}}(\sigma(a),\sigma(b))\leq b-a$. For any sequence $r_k<\bar{r}$, $r_k\rightarrow \bar{r}$, $\{\sigma(r_k)\}$ is a Cauchy sequence in  a \textbf{compact}  subset  $\bar{B}(p,L)$. So $\lim\limits_{s\rightarrow \bar{r}} \sigma(s)$ must exist. Denote the limit by $q$. 
 
 Choose a local coordinate system $\{x^\alpha\}$ as in $\S{2.1}$ around $q$. Since $\pi(\sigma(s))\triangleq \gamma(s)$, $s<\bar{r}$,  satisfies the 2nd  equation of the  ODE (\ref{pge}) near $q$, we know  $\sigma(s)$ can be extended smoothly over $\bar{r}$, i.e., $\sigma$ is now defined on $[0,\bar{r}+\epsilon]$ for some  $\epsilon>0$, and $\pi(\sigma)\mid_{[\bar{r}-\epsilon,\bar{r}+\epsilon]}$ satisfies  (\ref{pge}) on the coordinate system $\{x^{\alpha}\}$.  By solving the $x^0$-coordinate function from the 1st equation of ODE (\ref{pge}) for $s\in [\bar{r}-\epsilon, \bar{r}+\epsilon]$, we get a $\hat{g}$-geodesic $\tilde{\gamma}:[\bar{r}-\epsilon,\bar{r}+\epsilon]$ lying in the coordinate system whose horizontal projection curve is $\sigma\mid_{[\bar{r}-\epsilon,\bar{r}+\epsilon]}$. Since $X$ is complete (see Lemma \ref{Xcom}), one can choose a suitable $t_0\in R$, such that $\Psi_{t_0}(\tilde{\gamma})$ coincides with $\gamma$ on $[\bar{r}-\epsilon,\bar{r})$. Now $\gamma\cup \Psi_{t_0}(\tilde{\gamma}\mid_{[\bar{r},\bar{r}+\epsilon]}) $ will be a smooth $\hat{g}$- geodesic which is an extension of $\bar{\gamma}$. This is a contradiction with the definition of $\bar{r}$.    
 The proof is complete. 
 \end{pf}
 
\begin{thm} \label{t3.3}
Let $(M,g_{M})$ be  a geodesically complete static Einstein spacetime of dimension $n+1$, i.e. $Ric(g_{M})=\lambda g_{M}$.  Then $\lambda\leq 0$.
\end{thm}
\begin{pf}  Suppose $\lambda>0$.   We know $Ric(\hat{g})= \lambda \hat{g}$ by (\ref{hatric}).  On the other hand,  $(M,\hat{g})$ is complete by Theorem \ref{complete}.  This implies that $M$ is compact by Bonnet-Myers theorem.   From (\ref{lf}) and (\ref{equ}), we have $\hat{\triangle} \log u=-\lambda$ on $M$.  At the minimum point of $\log u$, we find $-\lambda\geq 0$, which is a contradiction with $\lambda>0$. 
\end{pf}

\section{Static solutions}
In this section, we will handle static spacetimes. We will first   derive a local gradient estimate on the norm of the Killing field.  The idea comes from  Yau's gradient estimate of  harmonic functions  on Riemannian manifolds (see \cite{Y} or $\S{1.3}$ in \cite{SY}).   
\subsection{Static vacuum solutions}

In the following Theorem \ref{staticgra},   we assume $(M,g_{M},X)$ is a  static Einstein spacetime of dimension $n+1$ with a timelike Killing field $X$ whose orthogonal complement is integrable,  and  $Ric(g_{M})=\lambda g_{M}$.

Let $\hat{g}$ be the  Riemannian metric  defined by  (\ref{hat}).   Note that  $u={[-g_{M}(X,X)]}^{\frac{1}{2}}$ is a time-independent function, we have $|\hat{\nabla} \log u|^2_{\hat{g}}=|\nabla \log u|^2_{g}$ by (\ref{inverse}).
 \begin{thm}\label{staticgra} Let $\hat{B}(x_0, a)$ be a $\hat{g}$-geodesic ball centered at $x_0$ of radius $a>0$ with compact closure  in $M$. Then there is a universal constant $C$ such that 
 \begin{equation} \label{gs}
\begin{split}
& \sup_{x\in \hat{B}(x_0,\frac{a}{2})}\mid\hat{\nabla}\log u\mid_{\hat{g}}\leq C(\sqrt{n}a^{-1}+\sqrt{\max\{-\lambda,0\}}).
\end{split}
\end{equation}
 \end{thm}
 \begin{pf}   
   By (\ref{covf}), we know 
 \begin{equation}\label{3.11}
 \begin{split}
 & \hat{\nabla}^2 \log u(e_0,e_0)=| \nabla u|^2\\
 & \hat{\nabla}^2\log u(e_0,e_i)=0\\
 &\hat{\nabla}^2\log u(e_i,e_j)= \nabla_{i}\nabla_j \log u,
 \end{split}
 \end{equation}
 on a local coordinate system $\{x^{\alpha}\}$ in $\S{2.1}$, where $e_0=\frac{\partial}{\partial t}$, $e_i=\frac{\partial}{\partial x^i}-\theta_i\frac{\partial}{\partial t},$ $i=1,2,\cdots, n$.  
 This implies $\hat{\triangle} \log u=-\lambda$.  
 
 By Corollary \ref{c3.1},  we know $\hat{g}$ is also an Einstein metric and  $\hat{R}ic=\lambda \hat{g}$.   
 
 By Bochner formula, we have 
  \begin{equation}\label{3.12}
 \begin{split}
 & \hat{\triangle} |\hat{\nabla} \log u|^2=2|\hat{\nabla}^2 \log u|^2_{\hat{g}}+2\lambda|\hat{\nabla} \log u|^2.
  \end{split}
 \end{equation}
 
 From (\ref{3.11}), we get 
 \begin{equation}\label{3.13}
 \begin{split}
 & |\hat{\nabla}^2 \log u|^2_{\hat{g}}=| \nabla^{2}_{ij}\log u|^2_{g}+|\nabla \log u|_{g}^4,
  \end{split}
 \end{equation}
 and 
  \begin{equation}\label{3.14}
 \begin{split}
 & \hat{\triangle} |\nabla \log u|^2=2|{\nabla}^2 \log u|^2_{{g}}+2|\nabla \log u|^4+2\lambda|\nabla \log u|^2. 
  \end{split}
 \end{equation}
 Let $\rho$ be the $\hat{g}$-distance function  centered at $x_0$.  Let $\psi:R_{+}\rightarrow R$ be a smooth nonnegative decreasing cutoff function such that $\psi=1$ on $[0,\frac{1}{2}]$, $\psi=0$ outside $[0,1]$, and  $|\psi^{''}|+\frac{(\psi')^2}{\psi}\leq C\sqrt{\psi}$. 
 
 We consider the nonnegative function $f=\psi(\frac{\rho}{a}) |\nabla \log u|^2$ on $\hat{B}(x_0,a)$. Suppose $f$ achieves its  maximum at some smooth point $x_1\in \hat{B}(x_0,a)$ of $\rho$. Then we have 
 $\hat{\triangle} f(x_1)\leq 0$ and $\hat{\nabla}f(x_1)=0$. Hence 
 \begin{equation}\label{3.15}
 \begin{split}
 & 0\ge \hat{\triangle} f(x_1)\ge 
 2\psi|{\nabla}^2 \log u|^2_{{g}}+2\psi |\nabla \log u|^4+2\lambda\psi |\nabla \log u|^2\\
 &\ \ \ -\frac{1}{a^2}(|\psi^{''}|+2\frac{(\psi^{'})^2}{\psi})|{\nabla} \log u|^2+a^{-1}\hat{\triangle}\rho \psi^{'} |\hat{\nabla} \log u|^2.
  \end{split}
 \end{equation}
 We first consider $\lambda\geq 0$ case. In this case, we have the  Laplacian   comparison theorem $\hat{\triangle}\rho\leq \frac{n}{\rho}$.  Hence $a^{-1}\hat{\triangle}\rho \psi^{'}\geq 2na^{-2}\psi^{'}$. Multiplying both sides of (\ref{3.15}) by $\psi$, we get 
 $2f(x_1)^2-{Cn}a^{-2}f(x_1)\leq 0$, which implies $f(x_1)\leq Cn{a^{-2}}$. In particular, we have 
 
  \begin{equation}\label{3.16}
 \sup_{x \in \hat{B}(x_0,\frac{a}{2})}|\hat{\nabla} \log u|_{\hat{g}}\leq C\sqrt{n}a^{-1}.
 \end{equation}
 
 If $x_1$ lies in the cut locus of $x_0$, by applying  a standard support function technique (see \cite{Y}, or Theorem 3.1 in \cite{SY}),  (\ref{3.15}) and (\ref{3.16}) still hold. 
 
 Now we assume $\lambda<0$. In this case, Laplacian  comparison theorem tells us $\hat{\triangle}\rho\leq \frac{n}{\rho}(1+\sqrt{\frac{|\lambda|}{n}}\rho)$ (see Corollary 1.2 in \cite{SY}). This gives  $a^{-1}\hat{\triangle}\rho \psi^{'}\geq (2na^{-2}+a^{-1}\sqrt{|\lambda|n})\psi^{'}$. Multiplying both sides of (\ref{3.15}) by $\psi$, we get $f(x_1)\leq C(|\lambda|+na^{-2})$, where $C$ is a universal constant (independent of $n$). The proof is complete. 
 \end{pf}

Proof of Theorem \ref{t1.2}.  If $(M,g_{M})$ is geodesically   complete, we know $(M,\hat{g})$ is complete from Lemma \ref{complete}. If $\lambda\geq 0$,  letting $a\rightarrow \infty$ in Theorem  \ref{staticgra},   one can prove  $u=const.$.  Now the 1-form $X^{\ast}$ dual to $X$ becomes closed. On the universal cover, $X^{\ast}=df$ must hold for some function $f$, the level set $\{f=const.\}$ is a global integrable submanifold of the horizontal distribution.  It is easy to see that  $\{f=const.\}$ is complete and Ricci flat.  The universal cover of $(M,g_{M})$ will be  isometric to $R\times \{f=const.\}$.  

The argument of Theorem \ref{staticgra} essentially provides a proof of the following: 

\begin{cor} Let $(M^{n}, g)$ be a Einstein Riemannian manifold  with a nowhere vanishing Killing field $X$ such that the orthogonal complement of $X$ is integrable and $Ric=\lambda g$.   Then for any metric ball $B(x_0,a)$ with compact closure in $M^{n}$, we have 
\begin{equation}\label{csg}
\begin{split}
\sup_{B(x_0,\frac{a}{2})}|\nabla \log |X| |\leq  C(\sqrt{n}a^{-1}+\sqrt{\max\{-\lambda,0\}}),
\end{split}
\end{equation}
where $C$ is a universal constant. Moreover, if $\lambda\geq 0$ and $(M^{n}, g)$ is complete, then 
 $|X|=const.$ and the universal cover of $(M^n,g)$ is isometric to $R\times N$, where $N$ is a complete Ricci flat Riemannian manifold.  
\end{cor}
 \begin{cor} Let $(N^n, g)$ be a Riemannian manifold, $u$ a smooth positive function on $N$ satisfying $R_{ij}=u^{-1}\nabla_{ij}u+\lambda g_{ij}$, $\triangle u=-\lambda u$.   Then for any metric ball $B(x_0,a)$ with compact closure in $N$, we have 
\begin{equation}\label{csg}
\begin{split}
\sup_{B(x_0,\frac{a}{2})}|\nabla \log u|\leq  C(\sqrt{n}a^{-1}+\sqrt{\max\{-\lambda,0\}}),
\end{split}
\end{equation}
where $C$ is a universal constant. Moreover, $u=const.$ and $(N^n,g)$ is Ricci flat if $\lambda\geq 0$ and $(N^n,g)$ is complete.  
\end{cor}
\begin{pf} Let $M=R\times N^n$ be a  manifold equipped with a static Riemannian metric $\hat{g}=u^2dt^2+g$. Now $X\triangleq \frac{\partial}{\partial t}$ is a Killing field. One can show that $Ric(\hat{g})=\lambda \hat{g}$ (see (\ref{Ric})). The same argument as in Theorem \ref{staticgra} will give (\ref{csg}).  Because the projection from $M$ to $N^n$ of any $\hat{g}$-Cauchy sequence on $M$ is also a $g$-Cauchy sequence on $N^n$, the completeness of $N^n$ will imply the completeness of $M$. The last assertion of the corollary holds. 
\end{pf}

\section{4-d stationary vacuum spacetimes}
\subsection{Preliminaries}
In this section,  we focus on the usual dimension of spacetime, i.e., dim $M=4$. Now fix a point $x_0\in M$,  let  $\{x^{\alpha}\}$ be a coordinate system  used in  $\S{2.1}$, which covers some open neighborhood $M'\subset M$ of $x_0$.  Let $\pi: M\rightarrow N$ be the projection from $M$ to the $X-$ orbit space $N$. Equip $N'=\pi(M')$ the horizontal Riemannian metric $g_{ij}$.    
Since $X^{\ast}=-u^2(dt+\theta)$, the Hodge dual of  $X^{\ast}\wedge dX^{\ast}$ (on $(M,\bar{g})$)  is $\pm u^3\ast d\theta$, where $\ast d\theta$ is the Hodge dual of $d\theta$ (on $(N',g)$).  Denote $\omega=u^{3}\ast d\theta$.     Note that $|\omega|^2=u^6|d\theta|^2=\frac{u^6}{2}|\Lambda|^2$, where the norm for a 2-form is taken by requiring $|e^1\wedge e^2|^2=1$ if $e^1,e^2$ are orthonormal.  Now $d^2\theta=0$ is equivalent to $d(u^{-3}\ast \omega)=0$, or $g^{ij}\nabla_{i}\omega_{j}=3\langle d\log u, \omega \rangle$.  

It should be noted that $\omega\otimes\omega$ is globally defined   on $M$ no matter whether  $M$ is orientable or not. 

Now we can rewrite equations (\ref{equ}) and (\ref{hatric}) as follows:  
 \begin{equation} \label{4d1}
\begin{split}
& R_{ij}={u^{-1}}{\nabla_{i}\nabla_j u}+\frac{1}{2}{u^{-4}}(\omega_i\omega_j-|\omega|^2 g_{ij})+\bar{R}ic(e_i,e_j)\\
& \triangle u =-\frac{1}{2}u^{-3}|\omega|^2+ u^{-1}\bar{R}ic(X,X) \\
& g^{kl} \nabla_k\omega_{l}=3g^{kl}\omega_{k}\nabla_{l}\log u\\
& (\ast d\omega)_j=\pm2u \bar{R}ic(X,e_j),
\end{split}
\end{equation} 
and 
 \begin{equation} \label{4d2}
\begin{split}
& \hat{R}ic({X,X})=u^{-2}| \omega|^2-\bar{R}ic(X,X) \\
& \hat{R}ic({X,e_j})=-\bar{R}ic(X,e_j)
\\& \hat{R}ic({e_i,e_j})=-u^{-4}(|\omega|^2g_{ij}-\omega_i\omega_j)+\bar{R}ic(e_i,e_j).
\end{split}
\end{equation} 
Recall that $\tilde{g}\triangleq u^2g$ in section 2.3, and we have (see (\ref{Rictil}) (\ref{equtil}) (\ref{4d1})): 
\begin{equation} \label{Rictil4}
\begin{split}
& \tilde{R}_{ij}=\frac{1}{2}u^{-4}\omega_{i}\omega_{j}+2\frac{u_iu_j}{u^2}+\bar{R}ic(e_i,e_j)-{u^{-2}\bar{R}ic(X,X)}g_{ij}\\
& u^{2}\tilde{\triangle}\log  u =-\frac{1}{2}u^{-4}| \omega|^2+ u^{-2}\bar{R}ic(X,X) \\
& \tilde{g}^{kl}\tilde{\nabla}_k\omega_{l}=4\tilde{g}^{kl}\omega_{k}\nabla_{l} \log u\\
&  (\ast_{\tilde{g}} d\omega)(\frac{\partial}{\partial x^j})=\pm2\bar{R}ic(X,e_j).
\end{split}
\end{equation}
By taking trace on the first equation of (\ref{Rictil4}), we find that   the scalar curvature $\tilde{R}$ of $\tilde{g}$ satisfies 
\begin{equation} \label{scalartil}
\begin{split}
& u^2\tilde{R}=\frac{1}{2}u^{-4}|\omega|^2+2u^{-2}{|\nabla u|^2}+\bar{R}-2{u^{-2}\bar{R}ic(X,X)},
\end{split}
\end{equation} 
where $\bar{R}$ is the scalar curvature of $(M,\bar{g})$.
\subsection{A map $\Phi$}
Throughout this subsection, we assume the following condition holds: 
 \begin{equation} \begin{split}\label{C}  \bar{R}ic(X, Y)=0 \ \ \text{whenever}  \ \   \bar{g}(X,Y)=0,\end{split} \end{equation}
 where $X$ is the timelike Killing field. In general, condition (\ref{C}) does not hold, while it holds when $(M,\bar{g})$ is Einstein, i.e., $Ric(\bar{g})=c\bar{g}$.

 When condition (\ref{C})  holds, from the last equation of (\ref{4d1}), we know  $\omega=d\psi$ holds locally for some function $\psi$ by Poincare lemma. In this case, we have 
 \begin{equation} \label{psi1}
\begin{split}
&\hat{\triangle} \psi=4\langle\nabla \psi, \nabla \log u\rangle_{\hat{g}},  \\    & \tilde{\triangle} \psi=4\langle\nabla \psi, \nabla \log u\rangle_{\tilde{g}}.
\end{split}
\end{equation}

Let $g_{-1}=y^{-2}{(dx^2+dy^2)}$ be the hyperbolic metric (sectional curvature $\equiv-1$) on  Poincare upper half plane $H=\{(x,y): x\in R, y>0\}$.  We define a map $\Phi: M'\rightarrow H$ by  $\Phi=(\psi,u^2)=(x,y)$.  Because $\psi$ and $u$ are time-independent, $\Phi$ is also a map from $N'=\pi(M')$ to $H$. 

\begin{lem} The map  $\Phi$  satisfies \\ 
\ \ \ i)  $\Phi^{\ast} g_{-1}=u^{-4}\omega\otimes\omega+4d\log u\otimes d\log u$; \\
\ \ \ ii)  $\hat{\triangle} \Phi=u^{2}\tilde{\triangle}\Phi=2\bar{R}ic(X,X)\frac{\partial}{\partial y}$, where $\hat{\triangle} \Phi$( or $\tilde{\triangle}\Phi$) is the harmonic map Laplacian between  two Riemannian manifolds $(M', \hat{g})$(or $(N',\tilde{g})$) and $(H,g_{-1})$. 
\end{lem}
\begin{pf}
We can calculate the Christoffel symbols  $\Gamma^{a}_{bc}$ of $g_{-1}$ as follows: 
\begin{equation} \label{ChH}
\begin{split}
& {\Gamma}_{11}^{1}={\Gamma}_{22}^{1}=\Gamma_{12}^{2}=0\\
& -\Gamma_{12}^{1}=-\Gamma_{22}^{2}=\Gamma^{2}_{11}={y^{-1}}.
\end{split}
\end{equation}

From the definition of the Hessian of a  map from $(M,\hat{g})$ to $(H,g_{-1})$, we have 
 \begin{equation} \label{hessPsi}
\begin{split}
& (\hat{\nabla}_{\alpha\beta} \Phi)^{a}= \hat{\nabla}_{\alpha\beta} \Phi ^{a}+\Gamma^{a}_{bc} \nabla_{\alpha}\Phi^{b}\nabla_{\beta}\Phi^{c}.
\end{split}
\end{equation}
Combining with (\ref{ChH}), it follows 
\begin{equation} \label{5.4}
\begin{split}
& (\hat{\nabla}_{\alpha\beta} \Phi)^{1}=  \hat{\nabla}_{\alpha\beta} \psi-\psi_{\alpha}(\log u^2)_{\beta}-\psi_{\beta}(\log u^2)_{\alpha}\\
& (\hat{\nabla}_{\alpha\beta} \Phi)^{2}= 2u^2 \hat{\nabla}_{\alpha\beta}\log u+ u^{-2}\psi_{\alpha}\psi_{\beta}.
\end{split}
\end{equation}

Taking traces  on (\ref{5.4}) with respect to $\hat{g}$,  we get 
\begin{equation}\label{laplacepsi}
\begin{split}
& (\hat{\triangle}\Phi)^{1}=  \hat{\triangle}\psi-4\langle \nabla \psi, \nabla \log u\rangle_{{g}}=0\\
& (\hat{\triangle} \Phi)^{2}= 2u^2 \hat{\triangle}\log u+ u^{-2}|\omega|^2_{{g}}=2\bar{R}ic(X,X),
\end{split}
\end{equation}
where we have used (\ref{psi1}) (\ref{4d1}) and (\ref{lf})). 
\end{pf}
\begin{cor} When $(M,\bar{g})$ is Ricci flat, the map $\Phi=(x,y)=( \psi, u^2)$ is a harmonic map from $(M', \hat{g})$ (or $(N, \tilde{g})$) to $(H,g_{-1})$. 
\end{cor}
Now we can apply the standard Bochner formula
\begin{equation}\label{Boc}
\begin{split}
 \hat{\triangle} e(\Phi)=& 2\langle \hat{\nabla} \Phi, \hat{\nabla} \hat{\triangle} \Phi \rangle+2|\hat{\nabla}_{\alpha\beta}\Phi|^2+2\langle\hat{R}ic, \Phi^{\ast}g_{-1}\rangle_{\hat{g}}\\& -2R_{abcd}\Phi^{a}_{\alpha}\Phi^{b}_{\beta}\Phi^{c}_{\gamma}\Phi^{d}_{\delta}\hat{g}^{\beta\delta} \hat{g}^{\alpha\gamma}\\
\end{split}
\end{equation}
where \begin{equation} \label{e}
\begin{split}
e(\Phi)=& {\hat{g}}^{\alpha\beta}(\Phi^{\ast}g_{-1})_{\alpha\beta}=u^{-4}|\omega|^2+4|\nabla\log u|^2\\ =& 2u^2\tilde{R}-(\bar{R}-2u^{-2}\bar{R}ic(X,X)).\end{split}
\end{equation}
Here  we have used (\ref{scalartil}).

Now we compute the term $I_1=\langle \hat{\nabla} \Phi, \hat{\nabla} \hat{\triangle} \Phi \rangle$ first.  

By Lemma 5.1 and (\ref{ChH}),  we have 
\begin{equation}\nonumber 
\begin{split}
\hat{\nabla}_{\alpha}\hat{\triangle} \Phi^a& =\frac{\partial}{\partial x^{\alpha}}(\hat{\triangle} \Phi^a)+\Gamma^{a}_{bc}\hat{\triangle}\Phi^b\frac{\partial \Phi^c}{\partial x^{\alpha}}\\
&=\frac{\partial}{\partial x^{\alpha}}(\hat{\triangle} \Phi^a)-2u^{-2}\bar{R}ic(X,X)\frac{\partial \Phi^a}{\partial x^{\alpha}},
\end{split}
\end{equation}
which implies 
\begin{equation}\nonumber
\begin{split}
\hat{\nabla}_{\alpha}\hat{\triangle} \Phi^1& =-2u^{-2}\bar{R}ic(X,X)\frac{\partial \psi}{\partial x^{\alpha}},\\
\hat{\nabla}_{\alpha}\hat{\triangle} \Phi^2& =\frac{\partial}{\partial x^{\alpha}}(2\bar{R}ic(X,X))-2u^{-2}\bar{R}ic(X,X) \frac{\partial u^2}{\partial x^{\alpha}}.
\end{split}
\end{equation}
Hence, 
\begin{equation}\label{l1}
\begin{split}
I_1=& \langle \hat{\nabla} \Phi, \hat{\nabla} \hat{\triangle} \Phi \rangle=-2u^{-2}\bar{R}ic(X,X)(u^{-4}|\omega|^2+4|\nabla \log u|^2)\\& +u^{-4}\langle \nabla(2\bar{R}ic(X,X)),\nabla u^2  \rangle.
\end{split}
\end{equation}
The 3rd term on the right hand side of (\ref{Boc}) can be computed as follows: 
\begin{equation}
\begin{split}\label{5.14}
\langle\hat{R}ic, \Phi^{\ast} g_{-1}\rangle_{\hat{g}}&
= [u^{-4}(\omega^{i}\omega^{j}-|\omega|^2g^{ij})+\bar{R}ic(e_k,e_l)g^{ik} g^{jl}] [u^{-4}\omega_{i}\omega_{j}+4\frac{u_iu_j}{u^2}]\\
&=-4u^{-6}(|\omega|^2 |\nabla u|^2-\langle \omega, \nabla u\rangle^2)\\& \ \ \ \ +g^{ik} g^{jl}(u^{-4}\omega_{i}\omega_{j}+4\frac{u_iu_j}{u^2})\bar{R}ic(e_k,e_l).
\end{split}
\end{equation}

Since $g_{-1}$ has constant sectional curvature $K\equiv -1$, we have 
\begin{equation}\begin{split}\label{5.15}
  & -R_{abcd}\Phi^{a}_{\alpha}\Phi^{b}_{\beta}\Phi^{c}_{\gamma}\Phi^{d}_{\delta}\hat{g}^{\alpha\gamma} \hat{g}^{\beta\delta}\\ & = [(\Phi^{\ast}g_{-1})_{\alpha\gamma}(\Phi^{\ast}g_{-1})_{\beta\delta}-(\Phi^{\ast}g_{-1})_{\alpha\delta}(\Phi^{\ast}g_{-1})_{\beta\gamma}]\hat{g}^{\alpha\gamma}\hat{g}^{\beta\delta}\\&=8u^{-6}(|\omega|^2 |\nabla u|^2-\langle \omega, \nabla u\rangle^2).
\end{split}
\end{equation}

Now we compute the term $|\hat{\nabla}_{\alpha\beta}\Phi|^2$ at a given point $(\bar{x},t_0)$. By definition,  $|\hat{\nabla}_{\alpha\beta}\Phi|^2=u^{-4}(|(\hat{\nabla}_{\alpha\beta}\Phi)^1|^2_{\hat{g}}+|(\hat{\nabla}_{\alpha\beta}\Phi)^2|^2_{\hat{g}})$.

 Let $\{x^{i}\}$ be a normal coordinate system  around  the fixed point $\bar{x}\in N$.     
Let $E_{0}=u^{-1}\frac{\partial}{\partial t}$, $E_{i}=\frac{\partial}{\partial x^i}-\theta_i\frac{\partial}{\partial t}$,  then $\{E_{\alpha}\}$ is an orthonormal basis of $\hat{g}$ at $(\bar{x},t_0)$, hence 
\begin{equation}
|(\hat{\nabla}_{\alpha\beta}\Phi)^a|^2_{\hat{g}}=\sum_{\alpha,\beta}[(\hat{\nabla}_{\alpha\beta}\Phi)^a]^2.
\end{equation}

From (\ref{covf}), we have 
\begin{equation} \label{5.17}
\begin{split}
&(\hat{\nabla}^2\log u)(E_0,E_0)=u^{-2}|\nabla u|_{g}^2\\
&(\hat{\nabla}^{2}\log u)(E_0,E_i)=-\frac{1}{2}\Lambda_{il}u_kg^{kl}=\pm\frac{u^{-3}}{2}\ast(\omega\wedge du)_{i}\\
&(\hat{\nabla}^{2}\log u)(E_i,E_j)= (\nabla^2\log u)(\frac{\partial}{\partial x^i},\frac{\partial}{\partial x^j}).
\end{split}
\end{equation} 
On the other hand, we have 
\begin{equation} \label{5.18}
\begin{split}
&(d\psi \otimes d\psi)(E_0,E_0)=(d\psi \otimes d\psi)(E_0,E_i)=0\\
&(d\psi \otimes d\psi)(E_i,E_j)= \psi_{i}\psi_{j}.\end{split}
\end{equation} 
From (\ref{5.4}), we get 
\begin{equation}\label{5.19}
u^{-4}|(\hat{\nabla}_{\alpha\beta}\Phi)^2|^2_{\hat{g}}=4|\nabla \log u|^4+2|u^{-2}\omega\wedge d\log u|^2+|2\nabla_{ij}\log u+u^{-4}\omega_i\omega_j|^2.
\end{equation}

To  compute the term \begin{equation}
\begin{split}\nonumber |(\hat{\nabla}_{\alpha\beta}\Phi)^1|^2_{\hat{g}}& =(\hat{\nabla}^2 \psi(E_0,E_0))^2+2\sum_{i}(\hat{\nabla}^2 \psi(E_0,E_i))^2\\ \ \ & \ \ +\sum_{i,j}[\hat{\nabla}^2\psi(E_i,E_j)-2\psi_i(\log u)_j-2\psi_j(\log u)_i]^2, \end{split}\end{equation}
we need (see (\ref{covf}))
\begin{equation} \label{5.20}
\begin{split}
& \hat{\nabla}^2\psi(E_0,E_0)=\langle \nabla \log u,\nabla \psi\rangle\\
&\hat{\nabla}^2\psi(E_0,E_i)=-\frac{1}{2}u\Lambda_{il}g^{lk}\nabla_k\psi=\pm\frac{u^{-2}}{2}\ast(\omega\wedge \omega)_{i}=0\\
&\hat{\nabla}^2\psi(E_i,E_j)=(\nabla^2\psi)(\frac{\partial}{\partial x^i},\frac{\partial}{\partial x^j}),
\end{split}
\end{equation}
hence 
\begin{equation}\label{5.21}
u^{-4}|(\hat{\nabla}_{\alpha\beta}\Phi)^1|^2_{\hat{g}}
=\langle d \log u,u^{-2}\omega\rangle^2
+u^{-4}|\nabla_i\omega_j-\omega_i(\log u^2)_j-\omega_j(\log u^2)_i|^2.
\end{equation}
Combining (\ref{5.21}) (\ref{5.19}) (\ref{5.15})(\ref{5.14})(\ref{l1}), we have 
\begin{equation}\label{Bo}
\begin{split}
\hat{\triangle}(\frac{1}{2}e(\Phi))&= \triangle (\frac{1}{2}e(\Phi))
+\langle\nabla \log u, \nabla(\frac{1}{2}e(\Phi))\rangle\\& =  4|\nabla \log u|^4 +\langle d \log u,u^{-2}\omega\rangle^2+|2\nabla_{ij}\log u +u^{-4}\omega_i\omega_j|^2\\& \ \ \ +u^{-4}|\nabla_i\omega_j-2\omega_i(\log u)_j-2\omega_j(\log u)_i|^2+6|u^{-2}\omega\wedge d\log u|^2\\
&\ \ \ +I_2+I_3,
 \end{split}
 \end{equation}
 where \begin{equation} \label{l23}
\begin{split}
\frac{1}{2}e(\Phi)&=\frac{1}{2}u^{-4}|\omega|^2+2|\nabla\log u|^2\\
I_2&=[\bar{R}ic(e_k,e_l)-2u^{-2}\bar{R}ic(X,X)g_{kl}]g^{ik} g^{jl} [u^{-4}\omega_{i}\omega_{j}+4\frac{u_iu_j}{u^2}]\\
l_3&=4u^{-2}\langle \nabla(\bar{R}ic(X,X)),\nabla \log u  \rangle.
\end{split}
\end{equation}

\subsection{A priori estimates} 
 The main result of this section is the following: 
\begin{thm}\label{N}Let $(M,g_M)$ be a   spacetime  of dimension $4$ with a timelike Killing field $X$.   Denote $\hat{g}=-\frac{2}{g_{M}(X,X)}X^{\ast} \otimes X^{\ast}+g_{M}$ the  Riemannian metric associated to $X$.  Let   $\hat{B}(x_0, a)$ be a $\hat{g}$-metric ball centered at $x_0$ of radius $a>0$ with compact  closure  in $M$, and assume \begin{equation}
\begin{split}
\label{cc}
\sup_{\hat{B}(x_0, a)}|Ric(g_{M})|_{\hat{g}}\leq a^{-2}.
\end{split}
\end{equation}
Then there is a universal constant $C$  such that 
 \begin{equation} \label{ce}
\begin{split}
\sup_{x\in \hat{B}(x_0,\frac{a}{2})}|\nabla \log u|^2+u^{-4}|\omega|^2 \leq \frac{C}{a^2},
\end{split}
\end{equation}
 where $u^2=-g_{M}(X,X)$.  Moreover, for any $p>1$, there is a constant $C_p>0$ depending only on $p$ such that  
\begin{equation}
\begin{split}
(\frac{1}{vol_{\hat{g}}(\hat{B}(x_0,\frac{a}{2}))}\int_{\hat{B}(x_0,\frac{a}{2})}|Rm(g_{M})|^p_{\hat{g}}dvol_{\hat{g}})^{\frac{1}{p}}\leq \frac{C_p}{a^2}
\end{split}
\end{equation}
 \end{thm}

\begin{pf} 
By scaling invariance, we may assume $a=1$.  The argument is divided into two cases: 

Case 1: $\partial \hat{B}(x_0,{a})\neq \phi$; 

 Case 2: $\partial \hat{B}(x_0,{a})=\phi$.

We treat Case 1 first.  Let $h(x)=2|\nabla \log u|^2(x)+\frac{1}{2}u^{-4}|\omega|^2(x)$, and  $f(x)=h(x)d_{\hat{g}}^2(x,\partial \hat{B}(x_0,1))$. Since $f$ is a nonnegative function on $\hat{B}(x_0,1)$, vanishes  on $\partial \hat{B}(x_0,1)$,  there is a point 
 $\bar{x}\in \hat{B}(x_0,1)$ such that $f(\bar{x})=\sup_{x\in \hat{B}(x_0,1)}f(x)$.  To prove the theorem,  it suffices to prove that there is a universal constant $C>0$ such that 
 $f(\bar{x})<C$. 
 We will argue by contradiction. Suppose there are  a sequence of spacetimes $(M_l,g_{M_l})$ and $\hat{g}_l-$balls $\hat{B}(x_l,1)\subset M_l$ with compact closure satisfying (\ref{cc}) with $a=1$,  but $f(\bar{x}_l)\rightarrow \infty$  as $l\rightarrow \infty$,   where 
\begin{equation}\nonumber
f(\bar{x}_l)=\sup_{x\in \hat{B}(x_l,1)} h_{l}(x)d^2_{\hat{g}_l}(x,\partial \hat{B}(x_l,1)),\ \  h_{l}(x)=2|\nabla \log u_l|^2+\frac{1}{2}u_l^{-4}|\omega_l|^2.
\end{equation} Now we will work on a fixed space $(M_l,\hat{g}_l)$. For simplicity, we drop the subscript $l$. 

For any fixed $0<\epsilon<1$,  and any $x\in \hat{B}(x_0,1)$ with $$d_{\hat{g}}(x,\bar{x})\leq \epsilon f(\bar{x})^{\frac{1}{2}} h^{-\frac{1}{2}}(\bar{x})=\epsilon d_{\hat{g}}(\bar{x},\partial \hat{B}(x_0,1)),$$
 we have 
 \begin{equation} \label{h}
 h(x)\leq \frac{1}{(1-\epsilon)^2} h(\bar{x}).
  \end{equation}
  
 Note that the function $f(x)$ is invariant under the scaling of the metric.  The metric  $\bar{g}=-u^2(dt+\theta)^2+g$ is  invariant under normalizations $u\rightarrow  {u(\bar{x})^{-1} u}$,  $t\rightarrow u(\bar{x}) t$, $\theta\rightarrow u(\bar{x})\theta$, $\omega\rightarrow u(\bar{x})^{-2} \omega$.  Therefore,  the equation (\ref{4d1}) remains invariant under such normalizations.  So without loss of generality, by scaling $u$ and the metric $\bar{g}$ by suitable positive constants,  we may assume  $u(\bar{x})=1$ and $h(\bar{x})=1$.  Now (\ref{h}) becomes 
 
  \begin{equation} \label{h1}
 h(x)\leq \frac{1}{(1-\epsilon)^2} \ \ \  \text{on}  \  B(\bar{x}, \epsilon f(\bar{x})^{\frac{1}{2}}). 
  \end{equation}
  From (\ref{h1}), we know $|\nabla \log u|\leq \frac{1}{\sqrt{2}}(1-\epsilon)^{-1}$ on $B(\bar{x}, \epsilon f(\bar{x})^{\frac{1}{2}})$. 
  
  Take $\epsilon=\frac{D}{\sqrt{f(\bar{x})}}$, where $D\geq 1$ is a fixed constant independent of $l$,  we have 
  \begin{equation} \label{u}
  e^{-\frac{D}{\sqrt{2}}(1-\frac{D}{\sqrt{f(\bar{x})}})^{-1}} 
  \leq u(x)\leq  e^{\frac{D}{\sqrt{2}}(1-\frac{D}{\sqrt{f(\bar{x})}})^{-1}} \ \ \  \text{on}  \  \hat{B}(\bar{x}, D).
  \end{equation}
  By (\ref{Rictil4}) and (\ref{u}), 
  it can be shown  that the sectional curvature of $\tilde{g}$ on  $\hat{B}(\bar{x}, D)$ satisfies \begin{equation}\label{Ktil}
 - \tilde{K}_{max}\leq  \tilde{K}(x)\leq e^{10D(1-\frac{D}{\sqrt{f(\bar{x})}})^{-1}}(1-Df(\bar{x})
  ^{-\frac{1}{2}})^{-2}  \triangleq \tilde{K}_{max}.
\end{equation}

Let $E=\mathcal{H}_{\bar{x}}\subset T_{\bar{x}}M $ be the orthogonal complement of $X$ at $\bar{x}$  equipped with the Euclidean metric induced from  $\bar{g}$ or $\hat{g}$. Let $\exp_{E}: E\rightarrow M$ be the restriction of the exponential map (of the conformal metric $u^2(\bar{g})$ or $u^2(\hat{g})$).  Clearly, $\exp_{E}$ is a smooth map. 

$\sl{Claim}$ 1:   There is  a universal constant $\delta_1>0$ such the exponential map $\exp_{E}$ is an immersion  from $B(0,\delta_1)\subset E$ to $\hat{B}(\bar{x}, 1)$. Moreover, the pull back $(0,2)$-tensor field  $\exp_{E}^{\ast} \tilde{g}$ is strictly positive definite everywhere on $B(0,\delta_1)$,  where  $\tilde{g}=u^2[\bar{g}-\frac{1}{g_{M}(X,X)}X^{\ast} \otimes X^{\ast}]$.

Let $v_1\in E$,  $0\neq v_2\in T_{v_1}E=E$, we will show  $(\exp_{E})_{\ast v_1}(v_2)\neq 0$ when $|v_1|$ is small.  Let ${\gamma}(s)=\exp_{E}(sv_1)$, $s\in [0, 1]$ be a horizontal geodesic w.r.t. $u^2(\bar{g})$  such that $\gamma(0)=\bar{x}$, $\gamma(1)=\exp_{E}(v_1)$. The variation $\epsilon\rightarrow \gamma_{\epsilon}(s)=\exp_{E}(s(v_1+\epsilon v_2)$ of horizontal geodesics  gives a  nontrivial Jacobi field $U=s(\exp_{E})_{\ast sv_1}(v_2)$ on $\gamma$ such that $U(0)=0, U(1)=(\exp_{E})_{\ast v_1}(v_2)$.  

Note that  one can always  construct  a contractible 3-dimensional smooth spacelike  immersed submanifold  $\sigma: \Sigma\rightarrow M$ so that  there is  a smooth map $\tilde{\gamma}:[0,1]\rightarrow \Sigma$ satisfying   $\gamma=\sigma\circ \tilde{\gamma}$. Actually, let $\{E_1(s),E_2(s),E_3(s)\}$, where $E_3(s)=\dot{\gamma}(s)$,  be a parallel and horizontal orthogonal  $u^2\bar{g}$-frame along  $\gamma$, let $\Sigma=[0,1]\times \{|w_1|^2+|w_2|^2<\epsilon^2\}\subset R^3$, then the map $\sigma: \Sigma\rightarrow M$, where $\sigma(s,w_1,w_2)\triangleq \exp_{E\gamma(s)}(w_1E_1(s)+w_2E_2(s))$,  will be an immersion when $\epsilon$ is sufficiently small.  By considering the integral curves of $X$ passing through  $\sigma(\Sigma)$ and setting the  affine parameters $t=0$ on $\Sigma$,  there is a small positive number $\delta>0$ such that the map  $F: \Sigma\times (-\delta, \delta)\rightarrow M$,  where $F(x,t)=\Psi_{t}(\sigma(x))$,   is also an immersion.  
Let $\pi: \Sigma\times (-\delta,\delta)\rightarrow \Sigma$, $\pi(x,t)=x$ be the natural projection map. Now we pull back the metric $\bar{g}$ by $F$ to $\Sigma\times (-\delta,\delta)$, and equip $\Sigma$ the horizontal metric(still denoted by $u^2 g$)  induced from  the metric $F^{\ast} (u^2 \bar{g})$. 

Note that we can lift the family of horizontal  geodesics $\epsilon \rightarrow \gamma_{\epsilon}$ on $M$  to a family of horizontal geodesics  $\epsilon \rightarrow \tilde{\gamma}_{\epsilon}$ on $\Sigma\times (-\delta, \delta) $ such that $F(\tilde{\gamma}_{\epsilon})={\gamma}_{\epsilon}$. Denote the variational vector field on $\tilde{\gamma}$ by $\tilde{U}$. Now $\pi(\tilde{\gamma}_{\epsilon})$
is a variation of geodesics on $(\Sigma,u^2g)$(see (\ref{pge})),  $\pi_{\ast} \tilde{U}$ is the variational  Jacobi field such that $\pi_{\ast} \tilde{U}(0)=0$. $\pi_{\ast} \tilde{U}$ is nontrivial since its derivative at $\tilde{\gamma}(0)$ with respect to $\dot{\tilde{\gamma}}(0)$  is  $\pi_{\ast}(v_2)\neq 0$.    By (\ref{Ktil}),  the conjugate radius of the exponential map of $( \Sigma, u^2 g)$  at $\tilde{\gamma}(0)$ is greater than ${\pi}\tilde{K}_{max}^{-\frac{1}{2}}$. So if $|v_1|< {\pi}\tilde{K}_{max}^{-\frac{1}{2}}$, we must have $\pi_{\ast} \tilde{U}(1)\neq 0$ and hence  $(\exp_{E}^{\ast}\tilde{g})_{v_1}(v_2)\neq 0$.  This finishes the proof of $\sl{Claim}$ 1. 

$\sl{Claim}$ 2:  The identity map from $E$ to itself is the exponential map at point $0$ of $(B(0,\delta_1), \exp_{E}^{\ast} \tilde{g})$.  

The  $\sl{Claim}$ 2 is clear by our construction. We  denote the metric $\exp_{E}^{\ast} \tilde{g}$ still by $\tilde{g}$. The point is that the injectivity radius of $\tilde{g}$ at $0$ is at least  $\delta_1$.  One  can also pull back the functions $u$,  the 1-form $\omega$  by the exponential map $\exp_{E}$ to $B(0,\delta_1)$. Since $B(0,\delta_1)$ is contractible, the function $\psi$ satisfying $d\psi=\omega$ can be globally defined  on $B(0,\delta_1)$. We denote these pulled back quantities  still by the same notations $u,\psi, \omega$.  From the first equation of  (\ref{Rictil4}) and our assumption, it is important to know that  the curvature of $\tilde{g}$ is bounded on $(B(0,\delta_1)$. By \cite{JK},  one can construct  a harmonic coordinate system $\{z^{i}\}$ of radius $2\delta_2>0$ around $0$ such that the estimate 
\begin{equation} \label{gc2}
\frac{1}{2}\delta_{ij}\leq\tilde{g}_{ij}\leq 2\delta_{ij}, 
|\tilde{g}_{ij}|_{C^{1,\alpha}}\leq \delta_2^{-1}
\end{equation}
 holds on  $\{|z|<2\delta_2\}$. 
 By the second and third equations of (\ref{Rictil4}), we know $\tilde{\triangle} u $ and $d\omega$ and $\delta \omega$ are uniformally bounded by our assumption. By elliptic regularity, $|u|_{C^{1,\alpha}}$ and $|\omega|_{C^\alpha}$ are uniformally bounded on a smaller ball.  By Arzela-Ascoli theorem, one can take a $C^{1,\alpha}$ and $W^{2,p}$ convergent subsequence for  $\tilde{g}$ and $u$, and  $C^{\alpha}$, $W^{1,p}$ convergent subsequence of $\omega$. We denote the limit by $\tilde{g}^{\infty}$, $u^{\infty}$, $\omega^{\infty}$.  Since $h$ only involves $\nabla u$,  $\omega$ and $g$, we know that the $C^{\alpha}$ norm of $h$ is uniformally bounded on $\{|z|<\frac{3}{2}\delta_2\}$(independent of $l$).  The limit $h^{\infty}$ is $C^{\alpha}$ and must satisfy
 \begin{equation}\label{maxp}
 h^{\infty}(x)\leq h^{\infty}(0)=1
 \end{equation}
 on $\{|z|<\frac{3}{2}\delta_2\}$. 
 Now we attempt  to show that the limit $(\tilde{g}^{\infty}, u^{\infty},\omega^{\infty})$ is actually smooth and  satisfies the vacuum Einstein equation.

 Recall that in the above harmonic coordinate system $\{z^{i}\}$, Ricci curvature $2\tilde{R}_{ij}=-\tilde{g}^{kl}\frac{\partial^2}{\partial z^k\partial z^l}\tilde{g}_{ij}+Q_{ij}(\partial \tilde{g},\tilde{g})$,
 where $Q$ is quadratic in $\partial \tilde{g}$, with   polynomial coefficients   in $\tilde{g}$, $\tilde{g}^{-1}$. 
 
For each scaled solution $g_l,u_l,\omega_l$, in the above harmonic coordinate system  $\{z^{i}\}$, multiplying the first equation of  (\ref{Rictil4}) by a  function $\xi \in W_0^{1,p}(B(0,\delta_2))$, $p>1$,  and integrating by parts, we get: 
 
 \begin{equation} \label{pre}
\begin{split}
& \int_{B(0,\delta_2)}\tilde{g}^{kl}\frac{\partial \tilde{g}_{ij}}{\partial z^k}\frac{\partial \xi}{\partial z^l}+[\partial_l\tilde{g}^{kl}\partial_k\tilde{g}_{ij}+Q(\partial \tilde{g}, \tilde{g})_{ij}]\xi dz^1dz^2dz^3\\ & =  2\int_{B(0,\delta_2)}(\frac{1}{2}u^{-4}\omega_{i}\omega_{j}+2\frac{u_iu_j}{u^2} )\xi dz^1dz^2dz^3+2I_4 
\end{split}
\end{equation}
 where 
 \begin{equation} \nonumber 
\begin{split}
& I_4=\int_{B(0,\delta_2)} ( \bar{R}ic(e_i,e_j)-{u^{-2}\bar{R}ic(X,X)}g_{ij})\xi dz^1dz^2dz^3.
\end{split}
\end{equation}
 Since the norm  of $\bar{R}ic(e_i,e_j)-{u^{-2}\bar{R}ic(X,X)}g_{ij}$ is bounded by $Ch(\bar{x}_l)^{-1}\rightarrow 0$ by our scaling, we know $I_4\rightarrow 0$ as $l\rightarrow \infty$.  Note that the $C^{\alpha}$-norms  of $\partial\tilde{g}$, $\omega$, $\partial u$ are uniformally bounded,   (\ref{pre}) must converge to 
 \begin{equation} \label{pos}
\begin{split}
& \int_{B(0,\delta_2)}(\tilde{g}^{\infty})^{kl}\frac{\partial}{\partial z^k}g^{\infty}_{ij}\frac{\partial}{\partial z^l}\xi \\& =\int_{B(0,\delta_2)} [-\partial_l(\tilde{g}^{\infty})^{kl}\partial_k\tilde{g}^{\infty}_{ij}-Q(\partial \tilde{g}^{\infty}, \tilde{g}^{\infty})_{ij}+(u^{\infty})^{-4}\omega^{\infty}_{i}\omega^{\infty}_{j}+4\frac{u^{\infty}_iu^{\infty}_j}{(u^{\infty})^2}]\xi.
\end{split}
\end{equation}
Because  $A_{ij}=-\partial_l(\tilde{g}^{\infty})^{kl}\partial_k\tilde{g}^{\infty}_{ij}-Q(\partial \tilde{g}^{\infty}, \tilde{g}^{\infty})_{ij}+(u^{\infty})^{-4}\omega^{\infty}_{i}\omega^{\infty}_{j}+4\frac{u^{\infty}_iu^{\infty}_j}{(u^{\infty})^2}\in W^{1,p},$  and coefficients $(\tilde{g}^{\infty})^{kl}\in W^{2,p}$, we know $\tilde{g}^{\infty}\in W^{3,p}$ by standard $L^p-$estimate for elliptic equations of divergence form. 

 Now we can apply the same technique to the rest  equations of  (\ref{Rictil4}) to obtain \begin{equation} \label{Rictil5}
\begin{split}
& (u^{\infty})^{2}\tilde{\triangle}\log  u^{\infty} =-\frac{1}{2}(u^{\infty})^{-4}| \omega^{\infty}|^2 \\
& (\tilde{g}^{\infty})^{kl}\tilde{\nabla}_k\omega^{\infty}_{l}=2(\tilde{g}^{\infty})^{kl}\omega^{\infty}_{k}\nabla_{l} \log u^{\infty}\\
&   d\omega^{\infty}=0
\end{split}
\end{equation}
where the above equations hold in the sense of integration by parts as in (\ref{pos}).
Since $\omega^{\infty}\in W^{1,p}$, $u^{\infty}\in W^{2,p}$, by applying $L^p$ estimates to the first equation of (\ref{Rictil5}), we know $u^{\infty}\in W^{3,p}$. This implies $\omega^{\infty}\in W^{3,p}$ by the 2nd and 3rd equations in (\ref{Rictil5}). Hence $A_{ij}\in W^{2,p}$, this gives $\tilde{g}^{\infty}_{ij}\in W^{4,p}$. Repeating this arguments, we find $\tilde{g}^{\infty}$, $u^{\infty}$, $\omega^{\infty}$ are actually smooth and satisfy the vacuum Einstein equations  on $\{|z|<\delta_2\}$. Since $d\omega^{\infty}=0$, from the calculations in section 5.2, we know  equation (\ref{Bo})  must hold for the limit $({g}^{\infty}, u^{\infty}, \omega^{\infty})$.   Moreover,  $I_2=0$ and $I_3=0$ hold in (\ref{Bo}) and (\ref{l23}).  That is to say,  we have 

 \begin{equation} \label{limit}
\begin{split}
& \triangle _{g^{\infty}}(h^{\infty})
+\langle\nabla \log u^{\infty}, \nabla(h^{\infty})\rangle\\ & =   4|\nabla \log u^{\infty}|^4+|\omega^{\infty}|^2|\nabla \log u^{\infty}|^2 +|2\nabla_{ij}\log u^{\infty} +(u^{\infty})^{-4}\omega^{\infty}_i\omega^{\infty}_j|^2 \\ & \ \ \  +(u^{\infty})^{-4}|\nabla_{i}\omega^{\infty}_j-2\omega^{\infty}_i(\log u^{\infty})_j-2\omega^{\infty}_j(\log u^{\infty})_i|^2\\& \ \ \ +5|(u^{\infty})^{-2}\omega^{\infty}\wedge d\log u^{\infty}|^2
\\& \  \geq 0
 \end{split}
 \end{equation}
 on $\{|z|<\delta_{3}\}$. 
 
Now we can apply the strong maximum principle on equation (\ref{limit}) since (\ref{maxp}) holds.  It follows that $h^{\infty}\equiv $ const., and the right hand side of (\ref{limit}) vanishes everywhere  on $\{|z|<\delta_{3}\}$. In particular, this implies $|\nabla \log u^{\infty}|^4\equiv 0$ and $|2\nabla_{ij}\log u^{\infty} +(u^{\infty})^{-4}\omega^{\infty}_i\omega^{\infty}_j|^2 $=0 on $\{|z|<\delta_{3}\}$, which give us  $u^{\infty}\equiv 1, \omega^{\infty}\equiv 0$. Hence $h^{\infty}\equiv 0$, which is a  contradiction with $h^{\infty}(0)=1$.  This proves Case 1. 

For Case 2, the maximum of $h(x)$ can be achieved at some point $\bar{x}$ by the compactness of $M$. The result can be proved by following  the same argument of Case 1.  The proof of the theorem is complete. 
 \end{pf}
 \vspace{3mm}
 
 From the proof of Theorem \ref{N},  if we integrate the vector field $X$ for a short time along   the image of the horizontal exponential map, one can obtain a local covering map which provides  a "good" local "coordinate system" (c.f.  (\ref{gc2})). 
 \begin{thm}\label{M} Under the assumptions of theorem \ref{N},  there is  a smooth non-degenerate  map $\varPsi: \{|z_0|^2+|z_1|^2+|z_2|^2+|z_3|^2<c^2 a^2\}\rightarrow \hat{B}(x_0, a)$, $\varPsi(0)=x_0$ such that $\varPsi^{\ast}g_{M}=\bar{g}_{\alpha\beta}dz^{\alpha}dz^{\beta}=-u^{2}(dz^0+\sum\theta_idz^i)^2+ g_{ij}dz^idz^j$ satisfies 
 \begin{equation}\label{lp5}
\begin{split}
& \frac{1}{1+C_0}<u<1+C_0, \ \  |\theta|<C_0\\
& (1+C_0)^{-1}\delta_{ij}<g_{ij}<(1+C_0)\delta_{ij}\\
&\frac{1}{a^4} \int_{\{|z|<ca\}}  (a|\partial \bar{g}|+a^2|\partial^2 \bar{g}|)^pdz<{C_p}
\end{split}
\end{equation}
 where $u$, $\theta$ and $g_{ij}$  are $z^0-$ independent,  $(z^1,z^2,z^3)$ are harmonic coordinates for $ g_{ij}dz^idz^j$,  $C_p$ are  constants  depending only on  nonnegative integers  $p\geq 0$. 
 \end{thm}

 \begin{thm} \label{A} Let $(M,g_{M})$ be  a  Einstein spacetime of dimension 4 with a timelike Killing field $X$,  $Ric(g_M)=\lambda g_{M}$, where $\lambda\geq 0$.   Let $\hat{B}(x_0, a)$ be  a $\hat{g}-$metric  ball in $M$ with compact closure. Then we have  \begin{equation}\label{5.38}
 \sup_{ \hat{B}(x_0,\frac{a}{2})} |\nabla \log u|^2+u^{-4}|\omega|^2\leq {C}a^{-2}, 
\end{equation}
for some universal constant $C$.
 \end{thm}
 
 \begin{pf} When $\lambda=0$, one can apply Theorem \ref{N} to derive (\ref{5.38}) since condition (\ref{cc}) holds trivially in this case.  We only need to handle $\lambda>0$ case.  By scaling invariance of the estimate, one can assume $a=1$.  We mimic the proof of Theorem \ref{N}.  
 
 We treat the case  $\partial \hat{B}(x_0,{1})\neq \phi$ first.   
 
 Let $h(x)=2|\nabla \log u|^2(x)+\frac{1}{2}u^{-4}|\omega|^2(x)+6\lambda$,  $f(x)=h(x)d_{\hat{g}}^2(x,\partial \hat{B}(x_0,1))$, and 
 $\bar{x}\in \hat{B}(x_0,1)$ such that $f(\bar{x})=\sup_{x\in \hat{B}(x_0,1)}f(x)$. To prove   $f(\bar{x})<C$ for some universal constant $C$,    we will argue by contradiction.  Suppose there are  a sequence of 4-Lorentzian manifolds $(M_l, \bar{g}_{l})$ satisfying ${R}ic(\bar{g}_l)=\lambda_{l}\bar{g}_l$ $(\lambda_l\geq 0)$  and a sequence of $\hat{g}_l$-balls $\hat{B}(x_l,1)\subset  M_{l}$  with compact closure such that $f(\bar{x}_l)\rightarrow \infty$  as $l\rightarrow \infty$,  where 
\begin{equation}\label{f}
\begin{split}
& f(\bar{x}_l)=\sup_{x\in \hat{B}(x_l,1)} h_{l}(x)d^2_{\hat{g}_l}(x,\partial \hat{B}(x_l,1))\\& h_{l}(x)=2|\nabla \log u_l|^2+\frac{1}{2}u_l^{-4}|\omega_l|^2+6\lambda_{l}.
\end{split}
\end{equation}
Scaling $u_l$ and $\bar{g}_l$ by $u_l(\bar{x}_l)^{-1}$ and $h_l(\bar{x}_l)$ respectively, one can assume  $u_l(\bar{x}_l)=1$, $h_l(\bar{x}_l)=1$.  We still use the same notations $u_l,$ $\omega_l$, $\bar{g}_l$, etc.,  to denote the corresponding scaled quantities.

Note that  the boundedness of $h_l$ implies that  the sectional curvature of $\tilde{g}$ is uniformally bounded on $B_{\hat{g}_l}(\bar{x}_l, 1)$.  As in Theorem \ref{N},  one can use  the horizontal exponential map (w.r.t. metric $u_l^2 \bar{g}_l$) to pull back $\tilde{g}_l$, $\omega_l$ and $u_l$ to horizontal tangent space.  Using  the harmonic coordinates $\{z^i\}$ on the horizontal tangent space and a boot strap argument as in Theorem \ref{N} , one can show that $\{u_l, \omega_l, \tilde{g}_l\}$ has  a subsequence converging to a smooth limit $(u^{\infty}, \omega^{\infty},\tilde{g}^{\infty})$.  Note that on (\ref{Bo}),   $I_2+I_3=4\lambda|\nabla \log u|^2+3\lambda u^{-4}|\omega|^2\geq 0$ for each $(u_l, \omega_l, \tilde{g}_l)$. So $I_2+I_3\geq 0$ still holds for the limit $\{u^{\infty}, \omega^{\infty},\tilde{g}^{\infty}\}$.   By  applying the strong maximum principle to equation  (\ref{Bo}) for the limit as in Theorem \ref{N}, we find $u^{\infty}\equiv 1$, $\omega^{\infty}=0$, $h^{\infty}\equiv 1$, and $\lambda^{\infty}=\frac{1}{6}$. From the second equation of (\ref{4d1}), we have $\triangle u^{\infty}=-\frac{1}{6}u^{\infty}$, which is a contradiction.  

If $\partial \hat{B}(x_0,{a})= \phi$,  $M$ will be compact. The maximum point of $h(x)$ can be achieved.    One can apply the strong maximum principle directly on (\ref{Bo}) to find a contradiction with $\lambda>0$ as in the preceding argument. 

\end{pf}
\begin{thm} \label{B}
 Let $(M,g_{M},X)$ be a  spacetime of dimension 4 with a timelike Killing field $X$ such that $Ric(g_M)=\lambda g_{M}$. Let $\hat{B}(x_0, a)$ be  a $\hat{g}$-metric  ball  with compact closure, where $0<a< \frac{1}{\sqrt{\max\{-\lambda,0\}}}$.  Then we have    \begin{equation} \label{ck2}
 \sup_{x\in \hat{B}(x_0,\frac{a}{2})} | Rm(g)|(x)+ | Rm(g_{M})|_{\hat{g}}(x)+|Rm(\hat{g})|_{\hat{g}}(x)\leq \frac{C_0}{a^{2}},
\end{equation}
and 
\begin{equation} \label{ck}
 \sup_{x\in \hat{B}(x_0,\frac{a}{2})} |\nabla_g^kRm(g)|_{g}(x)+ |\nabla_{g_M}^kRm(g_{M})|_{\hat{g}}(x)+|\nabla_{\hat{g}}^kRm(\hat{g})|_{\hat{g}}(x)\leq \frac{C_k}{a^{k+2}},
\end{equation}
where $Rm(g)$ is the Riemann curvature tensor of the horizontal metric $g$, $C_k,$ $k=0,1,2,\cdots$, are constants. 
\end{thm}
\begin{pf}
By scaling invariance, we can assume $a=1$, $\lambda\geq 0$ or $\lambda=-1$. By Theorems \ref{N} and \ref{A}, we know $|\nabla \log u|^2+u^{-4}|\omega|^{2}\leq  {C}$ on $\hat{B}(x_0,\frac{15}{16})$.  From this, we know the curvature $Rm(\tilde{g})$ of the conformal horizontal  metric $\tilde{g}=u^2g$ is bounded. Then we may apply the regularity argument  as in the proof of  Theorem \ref{N}  to prove 
\begin{equation} \label{ck1}
 \sup_{x\in \hat{B}(x_0, \frac{1}{2}+\frac{1}{4^{k+1}})} |\tilde{\nabla}^k u |_{\tilde{g}}(x)+ |\tilde{\nabla}^kRm(\tilde{g})|_{\tilde{g}}(x)+|\tilde{\nabla}^k \omega|_{\tilde{g}}(x)\leq C_k.
\end{equation}
Now (\ref{ck})  can be easily deduced from (\ref{ck1}) and (\ref{curvature}). 
\end{pf} 

Bing-Long Chen \\
 Department of Mathematics, \\
Sun Yat-sen University,\\
Guangzhou, P.R.China, 510275\\
Email: mcscbl@mail.sysu.edu.cn
\end{document}